\documentstyle[12pt]{article}
\title{The Chebyshev Exponent}
\author{Gene Ward Smith}
\date{}

\newcommand{\Q}{\mbox{$\bf Q$}}
\newcommand{\Z}{\mbox{$\bf Z$}}
\newcommand{\C}{\mbox{$\bf C$}}
\newcommand{\R}{\mbox{$\bf R$}}

\newcommand{\cp}{\scriptsize \copyright}
\newcommand{\eh}{\large \hat{e}}
\newcommand{\ec}{\large \check{e}}
\newcommand{\Sc}{\mbox{$\cal S$}}

\newcommand{\Uu}{\mbox{$\cal U$}}

\def\arccos{\mathop{\rm }}
\def\arccosh{\mathop{\rm arccosh}}

\def\Gal{\mbox{\rm Gal}}
\def\iso{\simeq}
\def\rtimes{\bowtie}

\def\d{\delta}
\def\e{\epsilon}

\newtheorem{theorem}{Theorem}
\newtheorem{proposition}[theorem]{Proposition}
\newtheorem{lemma}[theorem]{Lemma}
\newtheorem{corollary}[theorem]{Corollary}
\newtheorem{definition}{Definition}

\begin{document}

\maketitle

\begin{abstract}
The analogy between the nth power function and the nth Chebyshev polynomial
is pursued, leading to consideration of Chebyshev radicals as analogous to
ordinary radicals and Chebyshev exponents to ordinary exponents, and the
cosine and hyperbolic cosine as analogs of the exponential function. We
then discuss solving polynomial equations in Chebyshev radicals, and apply
this to the construction of unramified extensions of quadratic number fields.
\end{abstract}
 
\section{Basic Properties}

\subsection{Exponents}

There are many well-known families of polynomials $P_{n}(x)$
consisting of a polynomial of degree $n$ for each non-negative
integer $n$. Of all of these, the one which is at once the 
most important and the easiest to describe is the family of
powers, $P_{n}(x) = x^{n}$. This family may be defined by
means of a simple recurrence relationship, for we have
$P_{0}(x) = 1$, $P_{n}(x) = xP_{n-1}(x)$. It has also the 
remarkable property that functional composition is multiplicative
in the degree; that is
\[P_{n}(P_{m}(x)) = P_{nm}(x).\]

Suppose $E_{n}$ is a family of polynomials of the above sort,
so that $E_{n}$ is of degree $n$, and $E_{n}(E_{m}(x))=E_{nm}(x)$.
We may immediately make several observations. First, functional
composition commutes in such a family of polynomials, since
\[E_{n}(E_{m}(x))=E_{nm}(x) = E_{m}(E_{n}(x)).\]
Secondly, the family has a fixed point, namely, the constant value
of $E_{0}$, since
\[E_{n}(E_{0}(x)) = E_{0}(x).\]
Thirdly, we have
\[E_{1}(E_{1}(x)) = E_{1}(x),\]
and hence
\[E_{1}(x) = x.\]

Let $\ell(x) = px+q$ be a linear transformation. We may then
define a new family of polynomials by setting
\[\tilde{E}_{n} = \ell^{-1}(E_{n}(\ell(x)).\]
It is immediate that this family still has the functional
composition property,
\[\tilde{E}_{n}(\tilde{E}_{m}(x))=\tilde{E}_{nm}(x).\]
If we set $\ell(x) = (2x - b)/2a$, then
\[\tilde{E}_{n}(x) = aE_{n}(\frac{2x-b}{2a})+\frac{b}{2}.\]
This transformed family has the property that
if $E_{2}(x) = ax^{2}+bx+c$, then
\[\tilde{E}_{2}= x^{2} - \frac{b^{2}-4ac - 2b}{4}.\]

Let us therefore define an {\em exponent} to be a family of
polynomials $E_{n}(x)$, one for each degree $n$, such that
$E_{n}(E_{m}(x))=E_{nm}(x)$, normalized by the condition that
$E_{2}(x)$ is a monic polynomial with trace term zero; i.e.,
a polynomial of the form $x^{2}-c$. One can now ask, are there
in fact any exponents other than the familiar one, 
$P_{n}(x) = x^{n}$?  It turns out that there are;
in fact it follows from the work of Julia (\cite {Julia}) and Ritt
(\cite {Ritt}) that there are two and {\em only} two such exponents. 
(See also \cite {Riv} for a more elementary
treatment.)

The ``other'' exponent, of course, is the family of Chebyshev 
polynomials of the first kind.  Just as the family of power 
polynomials $P_{n}(x)$ can be defined by a recurrence relationship 
$P_{n}(x)=xP_{n-1}(x)$, so can these; we set $C_{0}(x) = 2$, $C_{1}(x) 
= x$, and \[C_{n}(x) = xC_{n-1}(x)-C_{n-2}(x).\]

The reader should note that other normalizations of 
this family are in use, the most common one in fact 
being $T_{0}(x) = 1$, $T_{1}(x) = x$, 
\[T_{n}(x) = 2xT_{n-1}(x)-T_{n-2}.\]
The relationship between these is simply
\[T_{n}(x) = \frac{C_{n}(2x)}{2}.\]
There are various advantages both conceptual and practical
to the normalization we have adopted, however; especially from
the point of view of an algebraist or number theorist.

\subsection{Chebyshev Powers}

\begin{definition}
For any non-negative integer $n$ and any $x$ in a ring $R$, we denote by
\[x^{\cp n}\]
the polynomial function defined by the
recurrence relation $x^{\cp 0} = 2$, $x^{\cp 1} = x$.
\begin{equation}
x^{\cp n} = x x^{\cp n-1} - x^{\cp n-2}. \label{eq:recur}
\end{equation}
For any negative integer $n$, we set $x^{\cp n} = x^{\cp -n}$.
\end{definition}

We will term this the {\em $n$-th Chebyshev power} of $x$.  It is, of 
course, a new notation for the $n$-th Chebyshev polynomial map, so that 
$x^{\cp n} = C_{n}(x)$.  The point is to emphasize the close analogy 
to the ordinary power $x^{n}$, and allow a more natural and suggestive 
notation for non-integral powers.  The reader who would like a way of 
pronouncing $x^{\cp n}$ is invited to pronounce it ``x cheby n''.

The special
symbol ``\copyright'' will be used to denote the Chebyshev analogue of
something expressed in terms of ordinary exponents. We will also
use it as a formal operator.
\begin{definition}
Let $R$ be a commutative ring and $R[x_{1}, \cdots,x_{n}]$ the
$R$-algebra of formal polynomials in $n$ indeterminates over $R$.
For $P \in R[x_{1},\cdots, x_{n}]$ we define $P \mapsto P^{\cp}$
to be the $R$-linear map on $R[x_{1},\cdots,x_{n}]$ which
takes each ordinary positive power $x_{i}^{k}$ in each monomial
term, and replaces it with the
corresponding Chebyshev power $x_{i}^{\cp k}$.
\end{definition}

Note that while the constant term of a polynomial 
$a_{0} x^{n}+ \cdots +a_{n}$ might be thought of as $a_{n}x^{0}$, 
we do {\em not} replace 
it by $a_{n}x^{\cp 0}=2 a_{n}$. So, for example, 
\[(x^{3}y^{2}+xy+1)^{\cp} = x^{\cp 3}y^{\cp 2}+xy+1.\]
Also, we must first expand the polynomial, so that for instance,
\[((x-1)(x^{2}+x+1))^{\cp}=(x^{3}-1)^{\cp} = x^{\cp 3}-1 = x^{3}-3x -1,\]
whereas
\[(x-1)^{\cp}(x^{2}+x+1)^{\cp}=(x-1)(x^{\cp 2}+x+1) = x^{3}-2x -1.\]

When an expression (for instance in $\Q(x)$) evaluates to a formal
polynomial, we will allow ourselves to operate on it with
``\copyright,'' so that for instance,
\[(\frac{x^{3}-1}{x-1})^{\cp} = x^{\cp 2} + x +1.\]
We will also make the usual identifications of formal polynomials
with polynomial functions when appropriate.

We may now derive a formula of fundamental importance, which in fact 
characterizes Chebyshev powers.

\begin{proposition}
Suppose $z \in R$ has an inverse element $z^{-1} \in R$ in a ring
$R$. 
For any integer $n$, we then have the relation
\begin{equation}
(z + z^{-1})^{\cp n} = z^{n} + z^{-n}. \label{eq:inv}
\end{equation}
\end{proposition}
Proof: By induction.
$\Box$

Note that in an algebraically closed field, such as $\C$, we may 
always find a $z$ satisfying $x = z+z^{-1}$ for any $x$, by solving 
$z^{2}-xz+1 = 0$. We can then get an expression \[x^{\cp 
n} = z^{n}+z^{-n}\] for $x^{\cp n}$ in terms of the roots of $z^2 - xz 
+ 1$, the characteristic polynomial for the recurrence relationship.

Equation (\ref{eq:inv}) also allows us to give another definition 
of the ``\copyright'' operator.
If $P(x) = a_{0}x^{n} + \cdots + a_{n}$ is a polynomial in
$R[x]$, where $R$ is commutative, then we set 
\[Q(z) = z^{n}(P(z) + P(z^{-1})-a_{0}).\]
We may then eliminate $z$ between two equations by taking the 
resultant of $Q(z)$ with $z^{2}-xz+1$; because of (\ref{eq:inv})
this will be $(P^{\cp})^{2}$.

\begin{proposition}
For positive integer exponent $n$,
\begin{equation}
x^{n} = \sum_{i \le n/2}{n \choose i}x^{\cp n-2i}. \label{eq:pow}
\end{equation}
\end{proposition}
Proof:
We can prove this by induction. We may also begin by
taking our ring be the complex numbers $\C$, and choosing $z$ such that
$x = z + z^{-1}$. Expanding $(z+z^{-1})^{n}$ by the binomial theorem and
collecting the terms $z^{i}$ and $z^{-i}$ together gives us the
result over $\C$. It then follows that the result is true as an
identity of formal polynomials over $\Z[X]$. If we send $X$ to
$x \in R$ and map $\Z$ to $R$, we have the result over any ring $R$.
$\Box$

\begin{proposition}
\label{primepow}
For any exponent $m = p^{n}$ which is a power of a prime $p$,
we have
\[x^{\cp m} \equiv x^{m} \pmod p.\]
\end{proposition}
Proof:
From the above, $x^{\cp p} \equiv x^{p} \pmod p$, which is to
say, modulo the ideal $pR$. By induction, 
\[x^{\cp n} = (x^{\cp n/p})^{\cp p} \equiv (x^{n/p})^{p}=x^{n} \pmod 
p.\]
$\Box$

\begin{definition}
We define a function $K(n,m)$ by
\[K(n, m) = {n \choose m} + {n-1 \choose m-1}.\]
\end{definition}

$K(n, m)$ can be expressed in a number of other ways. We have
\begin{eqnarray*}
K(n, m) = (1+\frac{n}{m}){n-1 \choose m-1} = 
\frac{n+m}{m}{n-1 \choose m-1} = \\
\frac{1}{m\!}(n+m)(n-1)(n-2) \cdots (n-m+1) = 
\frac{n+m}{n}{n \choose m}.
\end{eqnarray*}

For our purposes, the most important property of $K(n,m)$
is the following:
\begin{lemma}
\[K(n, m) = K(n-1, m) + K(n-1, m-1)\]
\end{lemma}
Proof: Immediate from the corresponding formula for
${n \choose m}$.

For positive $n$ we have $K(n,0) = 1$, $K(n, n) = 2$.
From this and the lemma, we see that we can define 
$K$ by means of a variant Pascal's triangle, with
1 along the left side, and 2 along the right. We have
\begin{center}
1 2 \\
1 3 2 \\
1 4 5 2 \\
1 5 9 7 2 \\
\end{center}
and so forth.

Our interest in $K(n, m)$ is a consequence of the following
proposition.
\begin{proposition}
For positive integer exponent $n$,
\begin{equation}
x^{\cp n} = \sum_{i \le n/2}(-1)^{i}K(n-i,i)x^{n-2i}. \label{eq:cheb}
\end{equation}
\end{proposition}
Proof:
By an easy induction.
$\Box$

\subsection{Chebyshev Radicals}

By far the most familiar method for solving polynomial equations
algebraically is the solution in radicals. 
A radical is an algebraic function defined by a choice of branch cut
for the polynomial $x^n - t$, giving $x$ as a function of $t$; so
that 
\[x = \sqrt[n]{t}.\]
This allows us to extend the definition of the exponent, originally 
defined for positive integers, to positive rational numbers; we are 
able to do this because the powers are exponents; that is, because
\[(x^n)^m = x^{nm}.\]

The $n$-th Chebyshev power map is in many respects analogous to the 
ordinary $n$-th power map.  If the ring $R$ is the complex numbers, 
then according to this analogy, an $n$-th Chebyshev root would be a 
solution of $x^{\cp n}-t$, and a Chebyshev radical would be a 
particular branch of this $n$-fold covering.  We will make a canonical 
choice of a branch which contains the fixed point $2^{\cp n}=2$.

\begin{definition}
For $n$ any non-zero integer,
we will denote by
\[\sqrt[\cp n]{x}\]
that branch of the function of a complex variable $x$ satisfying
\[(\sqrt[\cp n]{x})^{\cp n} = x\]
such that
$\sqrt[\cp n]{2} = 2$ and
with a branch cut along $(-\infty, -2)$, choosing
values along the cut to have positive imaginary part.

We will also denote this function by
\[x^{\cp {1 \over n}} = \sqrt[\cp n]{x}.\] 
\end{definition} 

\begin{proposition}
We have now the relation
\begin{equation}
\sqrt[\cp n]{\sqrt[\cp m]{x}} = \sqrt[\cp m]{\sqrt[\cp n]{x}} = 
\sqrt[\cp nm]{x}.
\end{equation}
\end{proposition}
Proof:
Since $(y^{\cp n})^{\cp m} = (y^{\cp m})^{\cp n} = y^{\cp nm}$, the
three expressions above all have the property that their $nm$-th
Chebyshev powers are equal to $x$. They also all have the same branch
cut $(-\infty, -2)$ and have positive imaginary part along the cut; 
hence they are equal. 
$\Box$

We now may now extend the definition to rational exponents, in analogy
to high school algebra.

\begin{definition}
For complex $x$ and $p/q$ any rational number, define
\[x^{\cp {p \over q}} = (\sqrt[\cp q]{x})^{\cp p}.\]
\end{definition}

\subsection{Chebyshev Exponentials and Logarithms}

In the case of ordinary powers, study of the powers of numbers close
to the fixed point for this family of exponents, $1^{n}=1$, is
particularly illuminating. If we look at the powers 
\[(1+\frac{1}{n})^{n}\]
we find that they tend to a particular number $e$ as $n$ grows
larger. Moreover, we find that if we define a function by
\[\exp z = \lim_{n \rightarrow \infty} (1+\frac{z}{n})^{n}\]
we produce something which may be regarded as the giving powers of
this number $e$, for any complex exponent. 

We may therefore define 
\[\exp x = e^x = \lim_{h \rightarrow 0} (1+ hx)^{1 \over h},\]
define a principal branch of the inverse function $\log$,
and discover that the relation
\[x^k = \exp (k \log x)\]
allows us to extend the definition of exponent to all complex exponents
$k$.

These two observations are connected. The exponential function 
satisfies the functional equation
\[\exp n x = (\exp x)^n,\]
and this in turn leads to the functional composition property of the 
exponent.

Perhaps the most illuminating way of extending the definition of the 
Chebyshev exponent to complex arguments is to employ a similar 
procedure to discover the Chebyshev analogs of the exponential and 
logarithmic functions.  We first define Chebyshev analogs to $e$ and 
$e^{-1}$, by looking at high powers of numbers near the fixed point 
for Chebyshev powers, $2^{\cp n}=2$.

\begin{theorem}
For a limit proceeding through rational values $h$, we have
\[\lim_{h \rightarrow 0}(2+h^{2})^{\cp{1 \over h}}=2 \cosh 1 = \eh.\]
\end{theorem}
Proof:
If we set
\[g = (h^{2} + \sqrt{h^{4}+4h^{2}})/2 = 
h+\frac{1}{2}h^{2}+\frac{1}{8}h^{3}+ \cdots\]
then
\[2+h^{2} = (1+g) + (1+g)^{-1}\]
and so
\[\lim_{h \rightarrow 0^{+}}(2+h^{2})^{\cp \frac{1}{h}}=\]
\[\lim_{h \rightarrow 0^{+}}(1+g+(1+g)^{-1})^{\cp \frac{1}{h}} =\]
\[\lim_{h \rightarrow 0^{+}}(1+h+O(h^2))^\frac{1}{h}+
(1+h+O(h^{2}))^{-\frac{1}{h}} = 
e+\frac{1}{e}.\]
$\Box$

\begin{proposition}
For $x$ any positive real number and $s$ any rational number, 
we have 
\[(x+x^{-1})^{\cp s} = x^{s}+x^{-s}.\]
\end{proposition}
Proof:
Let $s = p/q$, where $p$ and $q$ are integers, and let
$u = \sqrt[q]{x}$. Then
\[(u+u^{-1})^{\cp q} = u^{q}+u^{-q}=x+x^{-1}.\]
Since $u$ is positive and real, so is $u+u^{-1}$; which
therefore belongs to the correct branch and is the canonical
$n$-th Chebyshev root. Hence
\[u+u^{-1}=x^{1 \over q}+x^{-{1 \over q}}=x^{\cp {1 \over q}},\]
and so
\[x^{\cp s}=(x^{1 \over q}+x^{-{1 \over q}})^{\cp p}=
x^{p \over q}+x^{-{p \over q}} = x^{s}+s^{-s}.\]
$\Box$

We now have 
\[ \eh^{\cp s} = (e+e^{-1})^{\cp s} = e^{s}+ e^{-s} = 2 \cosh s \] 
for any rational number $s$. This suggests 
the following extended definition:
\begin{definition}
For any complex number $z$,
\[ \eh^{\cp z} = \exp_{\cp} z = 2 \cosh z .\]
\end{definition}

Despite the fact that it is essentially the hyperbolic cosine, we
introduce the notation $\exp_{\cp} z$ in order to underline the
formal similarities at work here. 

We have another limit analogous to the well-known
\[\lim_{h \rightarrow 0}(1-h)^{1 \over h} = e^{-1}\]
which is also of interest.
\begin{theorem}
Let $h$ approach 0 through rational values. Then
\[\lim_{h \rightarrow 0}(2-h^{2})^{\cp {1 \over h}}=2 \cos 1=\ec.\]
\end{theorem}
Proof:
Let $h = g/i$. Then
\[ \lim_{h \rightarrow 0}(2-h^{2})^{\cp {1 \over h}} = \]
\[ \lim_{g \rightarrow 0}(2+g^{2})^{\cp {i \over g}} = \]
\[ (\lim_{g \rightarrow 0}(2+g^{2})^{\cp {1 \over g}})^{\cp i} = \]
\[ \eh^{\cp i} = e^{i}+e^{-i} = 2 \cos 1 = \ec. \]

Now that we have identified $\eh^{\cp x} = 2 \cosh x$ and
$\ec^{\cp x} = 2 \cos x$ as the Chebyshev analogs of the
exponential function, we are led immediately to $\arccosh x/2$
and $\arccos x/2$ as the Chebyshev versions of the logarithm.
We make the following definition for the principal branch of the inverse
function to $\exp_{\cp} x$.
\begin{definition}
We define
\[ \log_{\cp}x = \arccosh \frac{x}{2}\]
for a branch such that if
\[w = \log_{\cp} = r+i \theta, \]
where $r = \Re w$ and $\theta = \Im w$, then
$r \ge 0$, $-\pi < \theta \le \pi$, and when $r=0$ then
$\theta \ge 0$.
\end{definition}

The reader may want to pronounce $\exp_{\cp}$ as
``cheby exp'' and $\log_{\cp}$ as ``cheby log''.

We may now extend our definition of Chebyshev powers
in a consistent way by setting
\begin{definition}
For any complex $x$ and $a$, we set
\[ a^{\cp x} = \eh^{\cp \log_{\cp}a x} = \exp_{\cp}(\log_{\cp}(a) x).\]
\end{definition}

In precalculus classes, students are told that for real and positive 
$a$, $a^{x}$ behaves in different ways depending on whether $a$ is 
less than, equal to, or greater than one.  In the same way, the 
behavior of $a^{\cp x}$ for $a \ge -2$ depends on whether it is less 
than, equal to, or greater than two; but in this case, the difference 
is more interesting, since for $a < 2$ we obtain a periodic function.

We now check that we have the properties that we want.
We note first that along the negative real axis from $-2$ to $-\infty$, 
$\log_{\cp}x$ is of the form $r + i \pi$, where $r$ is real 
and positive.  This gives us a positive imaginary part for the $n$-th 
Chebyshev root, as desired.

We now have
\begin{theorem}
For any complex $x$, $y$, and $a$,
\begin{equation}
(a^{\cp x})^{\cp y} = a^{\cp xy}. \label{eq:expo}
\end{equation}
\end{theorem}
Proof:
Since $f^{-1}(f(x)) = x$,
\[ f(y f^{-1}(f(x f^{-1}(a)))) = f(xy f^{-1}(a)). \]
Letting $f(x) = \eh^{\cp x}$, $f^{-1}(x) = \log_{\cp} x$, 
we have $a^{\cp x} = f(x f^{-1}(a))$, and hence the 
theorem upon substituting and simplifying.

One consequence is that
\[ \ec^{\cp x} = (\eh^{\cp i})^{\cp x} = \eh^{\cp ix}
= 2 \cosh ix = 2 \cos x, \]
as we suggested ought to be the case.

We may also now evaluate 
\[\lim_{h \rightarrow 0}(2+h^{2})^{\frac{1}{h}}\] 
by the same argument as before, with $h$ 
now any non-zero complex number.  If we substitute $h=x/n$, we then get 
\[\lim_{n \rightarrow \infty}(2+\frac{x^{2}}{n^{2}})^{\cp 
n/x}=\eh.\] 
Raising each side to the ``cheby x'' power, we get 
\[\lim_{n \rightarrow \infty}(2+\frac{x^{2}}{n^{2}})^{\cp 
n}=\eh^{x},\] 
and by restricting $n$ to integer values, we get an 
alternative definition of $\exp_{\cp} x$ 
analogous to the classical Cauchy definition of $\exp x$,
which simply requires us to 
take integral Chebyshev powers of numbers near to two.

We also have 
\begin{theorem}
\label{z}
For complex $z$ and $k$, let us define $z^{k}$ by means of the
principle branch of the logarithm, meaning the one such that
\[ -\pi < \Im(\log x) \le \pi.\]
Then except when
$z$ is a real number $-1 < z \le 0$, we have that
\[ (z + z^{-1})^{\cp k} = z^{k} + z^{-k}. \]
\end{theorem}
Proof:
Let $u = \log_{\cp}({z + z^{-1}})$. Then
\[ (z + z^{-1})^{\cp k} = \exp_{\cp} uk = e^{uk}+e^{-uk}. \]
If $z$ is not both real and negative, then one of $e^{uk}$ 
or $e^{-uk}$ will correspond to $z^{k}$, and the other to 
$z^{-k}$. If $z$ is real and less than $-1$, then $z = e^{u}$ and
$z^{k} = e^{uk}$, giving us the result we want. 
$\Box$

\subsection{The Derivative}

Throughout this section we will assume that the field of definition
is the complex numbers unless stated otherwise.

\begin{theorem}
Let $x$ and $k$ be any complex numbers, and let
$x=z+z^{-1}$, where $z$ is chosen so that it is not
the case that $-1 < z < 0$.
If $\Sc_{k}(x)$ is the function defined by 
\[\Sc_{k}(x) = \frac{z^{k}-z^{-k}}{z - z^{-1}},\] then
\[(x^{\cp k})' = k \Sc_{k}(x).\]
\end{theorem}
Proof:
If $x = z + z^{-1}$, then by Theorem \ref{z} we have
\[x^{\cp k} = z^{k} + z^{-k}.\]
Taking derivatives and applying the chain rule, we get
\[(x^{\cp k})' = \frac{kz^{k-1} - kz^{-k-1}}{1-z^{-2}} = k 
\frac{z^{k}-z^{-k}}{z - z^{-1}}.\]
$\Box$

Note that when $n$ is an integer, $\Sc_{n}(x)$ is a polynomial of 
degree $n-1$ in $x$, since it is a multiple of the derivative of a 
polynomial of degree $n$ in $x$.  Such polynomials are known as the 
Chebyshev polynomials of the second kind.
 
From the definition of $\Sc_{k}$ we immediately obtain
\[\Sc_{k}(\eh^{\theta}) = \frac{\sinh k \theta}{\sinh \theta},\]
\[\Sc_{k}(\ec^{\theta}) = \frac{\sin k \theta}{\sin \theta}.\]

We therefore have
\begin{equation}
\Sc_{k}(x) = \frac{\sinh (k \log_{\cp} x)}{\sinh \log_{\cp}x}
\end{equation}

\begin{theorem}
\label{S-recur}
\[\Sc_{k}(x) = x \Sc_{k-1}(x) - \Sc_{k-2}(x).\]
\end{theorem}
Proof:
\[\frac{(z+z^{-1})(z^{k-1}-z^{-k+1})-(z^{k-2}-z^{-k+2})}{z-z^{-1}} =
\frac{z^{k}-z^{-k}}{z-z^{-1}}.\]

For integer exponents $n$, we may consider $x^{\cp n}$ and 
$\Sc_{n}(x)$ to be elements of the function field $\C(x)$.  Consider 
any linear recurrence \[P_{n} = x P_{n-1} - P_{n-2},\] where 
$P_{n}(x)$ is a formal polynomial of degree $n$ in $\C(x)$.  Since 
this is a linear recurrence of degree two, the solutions comprise a 
vector space of dimension two over $\C(x)$.  Hence $x^{\cp n}$ and 
$\Sc_{n}(x)$ form a basis for the solutions of the linear recurrence.
Alternatively, we may write any such recurrence in terms of
$x^{\cp n}$ and $x^{\cp n-1}$, since these are also linearly
independent.

So, for instance, we have
\[x^{\cp n+1}=\frac{x x^{\cp n}+(x^{2}-4)\Sc_{n}(x)}{2},\]
or 
\[\Sc_{n}(x) = \frac{x x^{\cp n} - 2 x^{\cp n-1}}{x^{2}-4}\]
and so forth. Any of these may easily be derived if wanted by
using undetermined coefficients.

\subsection{Sum and Product Formulas}

For various purposes, it is helpful to have product formulas for 
our functions. We start with the Chebyshev analog of the product 
formula for powers, $x^{n}x^{m} = x^{nm}$. As before, arguments
$x$ and exponents $n$ or $m$ can be any complex number; results for
integral exponents can then be exported to any ring.

\begin{proposition}
\begin{equation}
x^{\cp n} x^{\cp m} = x^{\cp n+m} + x^{\cp n-m}. \label{eq:prod}
\end{equation}
\end{proposition}
Proof:
Writing $x$ as $x = z+z^{-1}$, we obtain
\[x^{\cp n} x^{\cp m} = (z^{n}+z^{-n})(z^{m}+z^{-m}) \]
\[=z^{n+m}+z^{-n-m}+z^{n-m}+z^{m-n}=x^{\cp n+m}+x^{\cp n-m}.\]

We can use (\ref{eq:prod}) to derive one of the most familiar 
properties of the Chebyshev polynomials, their orthogonality.

\begin{proposition}
Let $n, m$ be integers.
If $n \ne m$,
\[\int_{-2}^{2}x^{\cp n}x^{\cp m}\frac{dx}{\sqrt{4-x^{2}}} = 0,\]
whereas if $n \ne 0$
\[\int_{-2}^{2}(x^{\cp n})^{2}\frac{dx}{\sqrt{4-x^{2}}} = 2 \pi.\]
\end{proposition}
Proof:
Substituting $x=\ec^{\theta}$, we obtain
\[\int_{0}^{\pi}\ec^{\cp n \theta}\ec^{\cp m 
\theta}d\theta =
\int_{0}^{\pi}\ec^{\cp (n+m)\theta}+\int_{0}^{\pi}\ec^{\cp 
(n-m)\theta}d\theta.\]
This integral is 0 unless $n=m$, when it reduces to
\[\int_{0}^{\pi}2 = 2 \pi.\]

We also have a nice formula for the product of $\Sc_{n}$ with $x^{m}$.
\begin{proposition}
\begin{equation}
\Sc_{n}(x)x^{\cp m} = \Sc_{n+m}(x)+\Sc_{n-m}(x) \label{eq:pmix}
\end{equation}
\end{proposition}
Proof:
\[\frac{z^{n} - z^{-n}}{z-z^{-1}}(z^{m}+z^{-m})=
\frac{z^{n+m}-z^{-n-m} + z^{n-m}-z^{m-n}}{z-z^{-1}} = 
\Sc_{n+m}(x)+\Sc_{n-m}(x).\]

Finally, we have a formula for the product of $\Sc_{n}$ and $\Sc_{m}$.
\begin{proposition}
\begin{equation}
(x^{2}-4)\Sc_{n}(x)\Sc_{m}(x) = x^{\cp n+m} - x^{\cp n-m} \label{eq:ps}
\end{equation}
\end{proposition}
Proof:
\[\frac{(z^{n} - z^{-n})(z^{m}-z^{-m})}{z-z^{-1}} =
\frac{z^{n+m}+z^{-n-m}-z^{n-m}-z^{m-n}}{(z-z^{-1})^{2}} =
\frac{x^{\cp n+m} - x^{\cp n-m}}{x^{2}-4}.\]

\begin{corollary}
\begin{equation}
\Sc_{n}(x) = \frac{x^{\cp n+1}-x^{\cp n-1}}{x^{2}-4} 
\end{equation}
\end{corollary}
Proof:
From (\ref{eq:ps}), we have
\[(x^{2}-4)\Sc_{n}(x)\Sc_{1}(x) = x^{\cp n+1}-x^{\cp n-1}.\]
$\Box$

This or other expressions for $\Sc_{n}(x)$ in terms of Chebyshev
powers can be used to define it for any exponent.

We can use these product formulas to obtain formulas for the sum of
the exponents.

\begin{proposition}
\begin{equation}
x^{\cp n+m} = \frac{x^{\cp n}x^{\cp m}+(x^{2}-4)\Sc_{n}(x)\Sc_{m}(x)}{2}
\end{equation}
\begin{equation}
\Sc_{n+m}(x) = \frac{\Sc_{n}(x)x^{\cp m}+\Sc_{m}(x)x^{\cp n}}{2}
\end{equation}
\end{proposition}
Proof: From the product formulas.

These can be used to calculate high Chebyshev powers in any ring 
where two is invertible.

One way to regard these formulas is to see them as representing parts
of the powers of 
\[z = \frac{x + \sqrt{x^{2}-4}}{2},\]
so that
\[z^{n} = \frac{x^{\cp n} + \Sc_{n}(x)\sqrt{x^{2}-4}}{2}.\]

We can find a similar formula without any divisions.

\begin{proposition}
\begin{equation}
z^{n} = \Sc_{n}(x)z - \Sc_{n-1}(x)
\end{equation}
\begin{equation}
(x z - 2)z^{n-1} = (2 z-x)z^{n}=x^{\cp n}z - x^{\cp n-1}
\end{equation}
\end{proposition}
Proof: By induction.

Either of these propositions these can also be used to calculate 
Chebyshev powers by what is essentially the classic ``square and 
multiply'' algorithm for ordinary powers.  However, we can also avoid 
the use of anything other than Chebyshev powers themselves to 
calculate high-order Chebyshev powers.

\begin{lemma}
\[x^{\cp 2n} = (x^{\cp n})^{\cp 2}\] 
\[x^{\cp 2n+1} = x^{\cp n}x^{\cp n+1}-x\] 
\[x^{\cp 2n+2} = x x^{\cp 2n+1} - x^{\cp 2n}.\]
\end{lemma}
Proof: The first is the composition property of (Chebyshev) exponents, 
(\ref{eq:inv}).  The second follows from the product formula 
(\ref{eq:prod}).  The third is the recurrence relationship, 
(\ref{eq:recur}).

Using this, we may find any integral Chebyshev power by starting with 
the pair $(x^{\cp 1}, x^{\cp 2})$, and setting $l=n$.  At each stage, 
If $l$ is even, we replace $(x^{\cp m},x^{\cp m+1})$ with 
$(x^{\cp 2m}, x^{\cp 2m+1})$, and replace $l$ with $l/2$.  If 
$l$ is odd we replace $(x^{\cp m},x^{\cp m+1})$ with $(x^{\cp 
2m+1},x^{\cp 2m+2})$, and replace $l$ with $(l-1)/2$.  We 
continue in this way until $l=1$, when we arrive at $(x^{\cp n}, x^{\cp 
n+1})$.  If $n-1$ has fewer digits base two than does $n$, we may start 
instead by setting $l=n-1$ and finish with $(x^{\cp n-1}, x^{\cp n})$.

We now consider how, for integral exponents, we may express each of 
these polynomials in terms of the other.

\begin{proposition}
\begin{equation}
\Sc_{n+1}(x)-\Sc_{n-1}(x)=x^{\cp n} \label{eq:sn}
\end{equation}
\end{proposition}
Proof:
\[\Sc_{n+1}(x)-\Sc_{n-1}(x)=
\frac{z^{n+1}-z^{n-1}+z^{1-n}-z^{-1-n}}{z-z^{-1}}\]
\[= \frac{(z-z^{-1})z^{n} + (z-z^{-1})z^{-n}}{z-z^{-1}} = z^{n}+z^{-n} 
= x^{\cp n}.\]

\begin{proposition}
Let $n$ be integral.
For even $n$,
\begin{equation}
\Sc_{n}(x) = (\frac{x (x^{n}-1)}{x^{2}-1})^{\cp} = 
x^{\cp n-1} + x^{\cp n-3}+ \cdots + x, \label{eq:sev}
\end{equation}
while for odd $n$,
\begin{equation}
\Sc_{n}(x) = (\frac{x^{n+1}-1}{x^{2}-1})^{\cp} =
x^{\cp n-1}+x^{\cp n-3}+ \cdots + x^{\cp 2} + 1. \label{eq:sod}
\end{equation}
\end{proposition}
Proof:
On expressing the Chebyshev powers in terms of the $\Sc_{i}$ using 
(\ref{eq:sn}), we obtain a telescoping series.

We can use this to obtain a product formula for 
integral indices $n$ and $m$ of $\Sc_{n}\Sc_{m}$,
solely in terms of Chebyshev polynomials of the
second kind,  $\Sc_{i}$.
\begin{proposition}
Let $n$ and $m$ be positive integers. We then have the following:
\begin{equation}
\Sc_{n}\Sc_{m} = \sum_{i=1}^{m}\Sc_{n+m+1-2i}
\end{equation}
\end{proposition}
Proof: Suppose $m$ is even.  Expanding $S_{m}$ by means of 
(\ref{eq:sev}), we obtain 
\[\Sc_{n}(x)\Sc_{m}(x) = \Sc_{n}(x)\sum_{i=1}^{m/2}x^{\cp 2i-1}.\] 
Using 
(\ref{eq:pmix}), this becomes
\[\sum_{i=1}^{m/2}\Sc_{n+2i-1}+\Sc_{n-2i+1},\] 
which we may rearrange and reindex to obtain the theorem.

Similarly, if $m>1$ is odd, we may use (\ref{eq:sod}) and obtain
\[\Sc_{n}(x)\Sc_{m}(x)
=\Sc_{n}(x)(1+\sum_{i=1}^{(m-1)/2}x^{\cp 2i}).\]
Expanding this out gives
\[\Sc_{n}\Sc_{m}=\Sc_{n}+\sum_{i=1}^{(m-1)/2}\Sc_{n+2i}+\Sc_{n-2i},\]
and we can rearrange and reindex this to obtain the theorem. 
Since it is plainly true when $m=1$, we have a proof.
$\Box$

Finally, it is worth noting the connection of these to general 
second-order linear recurrences, and to Lucas and Fibonacci numbers 
and polynomials in particular.

Let 
\[U_{n} = PU_{n-1}-QU_{n-2},\]
\[V_{n} = PV_{n-1} - QV_{n-2},\]
be the ``fundamental'' and ``primordial'' Lucas sequences 
respectively, initialized by $U_{0}=0, U_{1}= 1$; $V_{0}=2,
V_{1}=P$. Then
\[U_{n}= Q^{\frac{n-1}{2}}\Sc_{n}(\frac{P}{\sqrt{Q}}),\]
\[V_{n}=Q^{\frac{n}{2}}(\frac{P}{\sqrt{Q}})^{\cp n}.\]

In particular, if $F_{0}=0,F_{1}=1, F_{n}=F_{n-1}+F_{n-2}$ are the 
Fibonacci numbers, and $L_{0}=2, L_{1}= 1, L_{n}= L_{n-1}+L_{n-2}$ are 
the Lucas numbers, then
\[\Sc_{n}(i) = i^{n-1}F_{n},\]
\[i^{\cp n} = i^{n} L_{n},\]
so that
\[F_{n}=i^{-n+1}\Sc_{n}(i) = i^{n-1}\Sc_{n}(-i),\]
\[L_{n}=i^{-n}i^{-\cp n} = i^{n}(-i)^{\cp n}.\]

We have in fact corresponding polynomials, where
\[F_{n}(x) = i^{n-1}\Sc_{n}(-ix)\]
are the Fibonacci polynomials, and
\[L_{n}(x) = i^{n}(-ix)^{\cp n}\]
are the Lucas polynomials.
These satisfy the recurrence 
relationship
\[P_{n}(x) = xP_{n-1}(x) + P_{n-2}(x).\]

The Lucas polynomials also satisfy the functional equation
\[L_{n}(z - z^{-1}) = z^{n} + (-z)^{-n},\]
and can be used in place of the Chebyshev polynomials to
solve algebraic equations.

\section{Algebraic and Analytic Properties}

\subsection{Some Factorizations}

For any univariate polynomial $p$ of positive degree, we have that
$x-y$ is a factor of $p(x)-p(y)$, since it is identically zero when 
$x=y$. When 
$p(x) = x^{n},$ we have the well-known expression for the difference
of two $n$-th powers:
\[x^{n} - y^{n} = (x-y)(x^{n-1}+x^{n-1}y + x^{n-2}y^{2} + \cdots + 
y^{n}).\] We seek the Chebyshev analogue to this; however we will
find it useful to look first at the factorization
for $\Sc_{n}(x)-\Sc_{n}(y)$.

\begin{proposition}
Let 
\[R_{n}(x,y) = \sum_{i=1}^{n-1}\Sc_{i}(x)\Sc_{n-i}(y).\]
Then
\begin{equation}
\Sc_{n}(x) - \Sc_{n}(y) = (x-y)R_{n}(x,y) \label{eq:sfac}
\end{equation}
\end{proposition}
Proof:
Since
\[x \Sc_{n}(x) = \Sc_{n+1}(x) + \Sc_{n-1}(x),\]
we have
\[(x-y)\Sc_{i}(x)\Sc_{j}(y) = \Sc_{i+1}(x)\Sc_{j}(y)-\Sc_{i}(x)\Sc_{j+1}(y)
+\Sc_{i-1}(x)\Sc_{j}(y)-\Sc_{i}(x)\Sc_{j-1}(y).\]
Hence
\[(x-y)R_{n}(x,y) = t_{1}-t_{2}+t_{3}-t_{4},\]
where
\begin{eqnarray*}
t_{1}=\sum_{i=1}^{n-1}\Sc_{i+1}(x)\Sc_{n-i}(y) = 
\Sc_{n}(x)+\sum_{i=2}^{n-1}\Sc_{i}(x)\Sc_{n+1-i}(y) \\
t_{2}=\sum_{i=1}^{n-1}\Sc_{i}(x)\Sc_{n+1-i}(y) = 
\Sc_{n}(y)+\sum_{i=2}^{n-1}\Sc_{i}(x)\Sc_{n+1-i}(y) \\
t_{3}=\sum_{i=1}^{n-1}\Sc_{i-1}(x)\Sc_{n-i}(y) = 
\sum_{i=1}^{n-2}\Sc_{i}(x)\Sc_{n-1-i}(y) \\
t_{4}=\sum_{i=1}^{n-1}\Sc_{i}(x)\Sc_{n-1-i}(y) =
\sum_{i=1}^{n-2}\Sc_{i}(x)\Sc_{n-1-i}(y)
\end{eqnarray*}
After cancelling the summation terms in the second expressions
given for the $t_{i}$, we are left with $\Sc_{n}(x)-\Sc_{n}(y)$.

\begin{theorem}
\begin{equation}
x^{\cp n} - y^{\cp n} = 
(x-y)(R_{n+1}(x,y)-R_{n-1}(x,y)) \label{eq:cfac}
\end{equation}
\end{theorem}
Proof:
By (\ref{eq:sn}), we have
\[x^{\cp n}-y^{\cp n}=\Sc_{n+1}(x)-\Sc_{n+1}(y)-\Sc_{n-1}(x)+\Sc_{n-1}(y).\]
Factoring this using (\ref{eq:sfac}) we obtain the proposition.

\begin{proposition}
\begin{equation}
x^{\cp n}-y^{\cp n} = 
(x-y)(\sum_{i=1}^{n-1}x^{\cp n-i}\Sc_{i}(y) + \Sc_{n}(y)) \label{eq:scfac}
\end{equation}
\end{proposition}
Proof:
Substituting $\Sc_{i+1}(x)-\Sc_{i-1}(x)$ for $x^{\cp i}$ in the above
and expanding, we obtain the result.

If $y=a$ is a constant, we can rewrite this as
\[x^{\cp n}-a^{\cp n} = (x-a)(\sum_{i=1}^{n}\Sc_{i}(a)x^{n-i})^{\cp}.\]
Special cases of this are of interest. If $a=-1$, then we have
\[\Sc_{3i}(-1)=0,\Sc_{3i+1}(-1)=1, \Sc_{3i+2}(-1) = -1.\]
From this we get
\[x^{\cp n}-(-1)^{\cp n} = 
(x+1)(x^{n-1}-x^{n-2}+x^{n-4}-x^{n-5}+ \cdots)^{\cp},\]
where the series on the right continues through all the 
non-negative exponents.

In a similar way, we can derive
\[x^{\cp n}-0^{\cp n} = x (x^{n-1}-x^{n-3}+x^{n-5} - \cdots)^{\cp},\]
\[x^{\cp n}-1^{\cp n} = (x-1)(x^{n-1}+x^{n-2}-x^{n-4}-x^{n-5}
+ \cdots)^{\cp},\]
and also
\[x^{\cp n}-(-2)^{\cp n}=(x+2)(x^{n-1}-2x^{n-2}+3x^{n-3}- 
\cdots)^{\cp},\]
\[x^{\cp n}-2^{\cp n} = x^{\cp n}-2 = (x-2)(x^{n-1}+2x^{n-2}
+3x^{n-3}+ \cdots)^{\cp}.\]

We also have formulas containing Fibonacci and Lucas numbers,
as for instance
\[x^{\cp n}-3^{\cp n}=x^{\cp n}-L_{2n} = (x-3)(x^{n-1}+F_{4}x^{n-2}
+F_{6}x^{n-3}+ \cdots)^{\cp}.\]

The factorization has some properties which correspond to the analogous 
factorization of the difference of ordinary powers. In particular,
we have the following.

\begin{theorem}
The factorization given by (\ref{eq:cfac}) is
irreducible over $\Q$ if and only if the exponent $n$ is prime.
\end{theorem}
Proof:
If $n = lm$ is not prime, then substituting $x^{\cp m}$ for $u$ and 
$y^{\cp m}$ for $v$ into the factorization of $u^{\cp l}-v^{\cp l}$ will 
further factorize $x^{\cp n} - y^{\cp n}$.

On the other hand, if $n$ is an odd prime, then $0^{\cp n} = 0$ shows
that 0 is a Chebyshev $n$-th power, and hence
\[x^{\cp n} - 0^{\cp n} = x (x^{\cp n-1} - K(n-1,1)x^{\cp n-3} + 
\cdots )\]
is the specialization obtained by setting $y=0$. 
Since $n$ is a prime, $x^{\cp n} = x^{n} \pmod n$, so all the 
coefficients in the second factor aside from the leading term are
divisible by $n$. The constant term is 
$(-1)^{(n-1)/2}K(\frac{n+1}{2}, \frac{n-1}{2})$
by (\ref{eq:cheb}).
Since 
\[K(m+1,m) = {m+1 \choose m}+{m \choose m-1} = 
{m+1 \choose 1} + {m \choose 1} = 2m+1,\]
we have 
\[K(\frac{n+1}{2},\frac{n-1}{2}) = 2 \frac{n-1}{2}+1 = n.\]
Hence the constant term is $(-1)^{\frac{n-1}{2}}n$, so that the second
factor is an Eisenstein polynomial, and hence irreducible. Since
a specialization is irreducible, the factorization in the theorem must be
irreducible also, and since the theorem is plainly true for the 
case $n=2$, we have a proof.

\subsection{Chebyshev Roots of Two}

For solving polynomial equations and for many other purposes, the
roots of unity $\zeta$ such that $\zeta^{n} = 1$ for some positive
integer $n$ are of fundamental importance. If the field in question
is the complex numbers, we have that $\zeta$ is a root of unity
if and only if
\[ \zeta = \exp (2 \pi i r)\]
for a rational number $r$. Here, if
$r = \frac{m}{n}$ we have $\zeta^{n} = \exp (2 \pi i m) = 1$.

The Chebyshev analog of this turns out to play a similar role when
solving equations using Chebyshev exponents. We first make
the following definition.
\begin{definition}
For any $a \in R$ in a ring $R$ and any integer $n$, we define
an $n$-th Chebyshev root of $a$ to be any element of $\mu \in R$ such 
that $\mu^{\cp n} = a$.
\end{definition}

In particular, a
{\em Chebyshev root of two} is an element $\mu$ such that
\[ \mu^{\cp n} = 2. \]

In close analogy with the situation for roots of unity, we have that when
the ring in question is the field of complex numbers,
\[ \mu = \ec^{\cp 2 \pi r} = 0^{\cp 4r} = \exp_{\cp}(2 \pi i r) \]
is a Chebyshev root of two if and only if $r$ is a rational number.
It follows that $\mu$ is real, and that $\mu = \zeta + \zeta^{-1}$, where
$\zeta$ is a root of unity.
If the {\em order} is the least positive $n$ such that
$\zeta^{n} = 1$ (for roots of unity) or $\mu^{\cp n} = 2$ (for 
Chebyshev roots of two) we also have that the order of $\mu$ and
the order of $\zeta$ are identical. If the order of $\zeta$ is
$n$, we say that $\zeta$ is a primitive $n$-th root of unity; in
the same way, we can call $\mu$ a primitive $n$-th Chebyshev root
of two if $n$ is the order of $\mu$.

If $x$ is an $n$-th Chebyshev root of two, then 
$x^{\cp n} - 2 = 0$.
If we factor this using (\ref{eq:scfac}), we obtain
\[x^{\cp n}-2 = (x-2)(\sum_{i=1}^{n}ix^{n-i})^{\cp}.\]

If $n>2$ this factorization is further reducible. We 
have already remarked that composite exponents such as
$x^{\cp 2m} - y^{\cp 2m}$ (for $m>1$) can be factored further, and
in fact we have the following.
\begin{proposition}
For any complex $x$ and $k$, we have
\begin{equation}
\label{eq:evefact}
x^{\cp 2k} - 2 = (x^{2}-4)\Sc_{k}^{2}(x).
\end{equation}
\end{proposition}
Proof:
Immediate from (\ref{eq:ps}).
$\Box$

For odd integer exponents, we have the following.
\begin{theorem}
For any positive integer $n$,
\begin{equation}
x^{\cp 2n+1} - 2 = (x-2)((\sum_{i=0}^{n}x^{n-i})^{\cp})^{2} 
= (x-2)(x^{\cp n} + x^{\cp n-1} + \cdots + x +1)^{2} \label{eq:cyc}
\end{equation}
\end{theorem}
Proof:
If $x$ is a $2n+1$-th Chebyshev root of two other than two,
and if $x = z+z^{-1}$,
then $z$ is a $2n+1$-th root of unity other than one, and so
\[z^{2n}+z^{2n-1}+ \cdots + 1 = 0.\]
Dividing by $z^{n}$ and collecting terms we obtain
\[z^{n} + z^{-n}+ z^{n-1}+z^{-n+1}+ \cdots + 1 = \]
\[x^{\cp n} + x^{\cp n-1} + \cdots + 1 = 0.\]
Since we obtain the root $x$ for each of $z$ and $z^{-1}$, it
appears with multiplicity two, and hence the factor is squared.
$\Box$

We can rewrite this factorization as the Chebyshev analog of 
a geometric series, for if $x \ne 2$ we have
\[ (1+ x + x^{\cp 2} + \cdots + x^{\cp n})^{2}=\frac{x^{\cp 
2n+1}-2}{x-2}.\]

As is so often the case, we can also express this factorization more
compactly using Chebyshev polynomials of the second kind, $\Sc_{n}$.
To this end and others, we introduce a new function.
\begin{definition}
For any complex $k$ and $x$ we define $\Uu_{k}(x)$ by
\[\Uu_{k}(x) = \Sc_{(k+1)/2}(x) + \Sc_{(k-1)/2}(x).\]
\end{definition}

\begin{proposition}
We have a recurrence relationship
\begin{equation}
\label{eq:u-recurr}
\Uu_{k+4}(x) = x \Uu_{k+2}(x) - \Uu_{k}(x)
\end{equation}
\end{proposition}
Proof:
From the recurrence relationship for $\Sc_{k}$.
$\Box$

\begin{corollary}
For positive integral $n$, $\Uu_{2n+1}(x)$ is given by the recurrence
relationship $\Uu_{1}(x) = 1$, $\Uu_{3}(x) = x+1$,
\[\Uu_{2n+3}(x) = x \Uu_{2n+1}(x) - \Uu_{2n-1}(x).\]
\end{corollary}

\begin{proposition}
For any complex $k$ and $x$, we have
\begin{equation}
\label{eq:sup}
\Uu_{k}(x) = \sqrt[\cp]{x}\Sc_{k/2}(x) = \sqrt{x+2}\Sc_{k/2}(x).
\end{equation}
\end{proposition}
Proof:
Apply (\ref{eq:pmix}) to the right-hand side.
$\Box$

We now have
\begin{proposition}
For any complex $x$ and $k$, we have
\begin{equation}
\label{eq:oddfact}
x^{\cp 2k+1}-2 = (x-2)\Uu_{2k+1}^{2}(x).
\end{equation}
\end{proposition}
Proof:
From (\ref{eq:evefact}), we obtain
\[x^{\cp 2k+1}-2 = (x^{2}-4)\Sc_{k+1/2}^{2}(x).\]
From (\ref{eq:sup}), we have
\[(x^{2}-4)\Sc_{k+1/2}^{2}(x) = (x+2)\Uu_{2k+1}(x),\]
and hence the theorem.
$\Box$

For integral $n$, we have a nice characterization of $\Uu_{2n+1}$
in terms of its roots.
\begin{proposition}
For positive integers $n$, we have
\[\Uu_{2n+1}(x) = \prod_{i=1}^{n}(x-\mu_{i});\]
where the $\mu_{i}$ are the $n$ $2n+1$-th Chebyshev roots of two other
than two.
\end{proposition}
Proof:
From (\ref{eq:oddfact}), we know a root of $\Uu_{2n+1}(x)=0$ must be
a $2n+1$-th Chebyshev root of two. From the recurrence relationship
we have that $\Uu_{2n+1}(2)=2n+1$, so two is not a root; the other 
$n$ roots are distinct from two and give the factors.
$\Box$

We may do something similar for $\Sc_{n}(x)$.
\begin{proposition}
If $n>1$ is integral, then
\[\Sc_{n}(x) = \prod_{i=1}^{n-1}(x-\mu_{i});\]
where the $\mu_{i}$ are the $n-1$ $2n$-th Chebyshev roots of two other
than $\pm 2$.
\end{proposition}
Proof:
From (\ref{eq:evefact}), we know a root of $\Sc_{n}(x)=0$
must be a $2n$-th root of two.
Since $\Sc_{n}(\pm 2) = \pm n$, all of the roots of 
$\Sc_{n}(x)$ must be distinct from $\pm 2$.
$\Box$

Finally, we may characterize $x^{\cp n}$ in the same way.
\begin{proposition}
For positive integers $n$,
\[x^{\cp n}= \prod_{i=1}^{n}(x-\mu_{i});\] 
where the $\mu_{i}$ are the $n$ $4n$-th Chebyshev roots
of two which are not also $2n$-th roots of two.
\end{proposition}
Proof:
If $\mu^{\cp n}=0$, then $\mu^{\cp 4n}= 0^{\cp 4} = 2$, and
hence any root of $x^{\cp n}=0$ must be a $4n$-th root of two.
However if $x^{\cp 2n}=2$, so that $x$ is also a $2n$-th root of
two, then $x^{\cp n} = \pm 2$. Hence $\mu$ must be a $4n$-th
root of two which is not a $2n$-th root of two; there are $n$ of
these, which is the degree of $x^{\cp n}$, so it must have all of them
as roots.
$\Box$

\begin{theorem}
For positive integral $n$, we have the factorizations
$$\Sc_{2n}(x) = \Sc_{n}(x) x^{\cp n}$$
$$\Sc_{2n+1}(x) = (-1)^{n}\Uu_{2n+1}(x)\Uu_{2n+1}(-x)$$
\end{theorem}
Proof:
The first formula is an immediate consequence of ({\ref{eq:pmix}).
The second follows from the fact that since
$x^{\cp 4n+2}-2 = (x^{\cp 2n+1}+2)(x^{\cp 2n+1}-2)$,
the roots of
$\Uu_{2n+1}(-x)$ are the $2n+1$-th roots of $-2$ other
than $-2$, which is to say the $4n+2$-th roots of two which are not
also $\pm 2$ or $2n+1$-th roots of two.
$\Box$

For Chebyshev powers, we have the following.
\begin{theorem}
For positive integral $n$, 
\begin{equation}
\label{oddu}
x^{\cp 2n+1} = (-1)^{n}x \Uu_{2n+1}(2-x^{2})
\end{equation}
\end{theorem}
Proof:
As before, the roots of $\Uu_{2n+1}(-x)$ are $2n+1$-th Chebyshev roots
of $-2$ other than $-2$; hence the roots of $\Uu_{2n+1}(-x^{\cp 2})$
are $2n+1$-th Chebyshev roots of $0$ other than $0$. Adding the
factor $x$ to get the root $x=0$ and multiplying by $(-1)^{n}$ to
ensure the polynomial is monic, we obtain the theorem.
$\Box$

For exponents which are powers of two, we have the following.
\begin{theorem}
\label{irrtwo}
If $n = 2^{m}$ is a power of two,
the polynomial $x^{\cp n}$ is irreducible over $\Q$.
\end{theorem}
Proof:
We have
\[x^{\cp n} \equiv x^{n} \pmod 2 \]
by Proposition \ref{primepow}. Hence each of the coefficients of
$x^{\cp n}$ other than the leading term is divisible by two.
Since $x^{\cp n} = (x^{\cp n/2})^{2}-2$, the constant
term is divisible by two only once, and so the polynomial is
Eisenstein, and hence irreducible.
$\Box$

Putting (\ref{oddu}) and (\ref{irrtwo}) together, we have the following.
\begin{corollary}
If $n = lm$, where $l$ is a power of two and $m$ is odd, then
\[x^{\cp n} = x^{\cp l}\Uu_{m}(-x^{\cp l}).\]
\end{corollary}

We therefore may conclude that the problem of factoring our various
polynomials reduces to the problem of factoring $\Uu_{2n+1}(x)$.

\begin{corollary}
The $2n+1$-th Chebyshev roots of two other than two are
units in the maximal real subfield of $2n+1$-th roots of
unity.
\end{corollary}

\begin{proposition}
If $n$ is not divisible by four, the $n$-th Chebyshev roots of unity
other than $\pm 2$ are units in the maximal real subfield $\Q(\mu)$
of the cyclotomic field $\Q(\zeta)$, where $\zeta$ and $\mu$ are
primitive $n$-th roots of unity and primitive Chebyshev $n$-th roots
of two, respectively.
\end{proposition}

The irreducible factors of $x^{n} - 1$ are the cyclotomic polynomials.  
The cyclotomic polynomial (of degree $\phi(n)$) whose roots are the 
primitive roots of unity of order $n$ is called the $n$-th cyclotomic 
polynomial, and is ordinarily denoted $\Phi_{n}$. If $\zeta$ is
a primitive $n$-th root of unity, then so is $\zeta^{-1}$. Hence
the coefficients of $\Phi_{n}$ for $n>1$ are palindromic; if
\[\Phi_{n}(x)=\sum_{i=0}^{\phi(n)}a_{i}x^{\phi(n)-i}\]
then $a_{i} = a_{\phi(n)-i}$. Since $\phi(n)$ is even when $n>2$, we may
make the following definition.

\begin{definition}
For $n>2$ we define the {\em $n$-th Chebyshev-cyclotomic polynomial}
to be
\[\Psi_{n}(x) = 
(\sum_{i=0}^{\phi(n)/2}a_{i}x^{\phi(n)/2-i})^{\cp},\]
where
\[\Phi_{n}(x) = \sum_{i=0}^{\phi(n)}a_{i}x^{\phi(n)-i}\]
is the $n$-th cyclotomic polynomial.
We also set $\Psi_{1}(x) = 1$.
\end{definition}

We then have
\begin{theorem}
For $n>2$,
$\Psi_{n}$ is an irreducible polynomial of degree $\phi(n)/2$
whose roots are the $\phi(n)/2$ primitive $n$-th Chebyshev roots
of two.
\end{theorem}
Proof:
Since $n>2$, $\phi(n)$ is even. If $x = z+z^{-1}$ and 
$\Phi_{n}(z) = 0$, we may divide by $z^{\phi(n)/2}$ and
collect terms, obtaining
\[\sum_{i=0}^{\phi(n)/2}a_{i}z^{\phi(n)/2-i}+a_{n-i}z^{i-\phi(n)/2}.\]
Since $a_{i}=a_{\phi(n)-i}$, this can be rewritten $\Psi_{n}(x)$.
If $i$ is relatively prime to $n$, $z^{i}$ is a Galois
conjugate of $z$; the roots of $\Phi_{n}(z)$ are the $\phi(n)$
Galois conjugates $z^{i}$ for $0<i<n$ prime to $n$, and so $\Phi_{n}$
is irreducible. In the same way, $x^{\cp i}$ for $i$ prime to
$n$ is a root of $\Psi_{n}$. If we take all $i$ prime to $n$
such that $0 < i < n/2$, we obtain a full
orbit of conjugates in number equal to the degree $\phi(n)/2$ of
$\Psi_{n}$, and hence $\Psi_{n}$ must be irreducible.
$\Box$

\begin{theorem}
For any  positive integer $n$, we may factor $\Uu_{2n+1}(x)$ by
\begin{equation}
\Uu_{2n+1}(x) = \prod_{d | 2n+1}\Psi_{d}(x)
\end{equation}
\end{theorem}
Proof:
Each side of the equation is a monic polynomial whose roots are
all of the $2n+1$-th Chebyshev roots of two other than two.
$\Box$

\subsection{Differential Equations}

The product formula (\ref{eq:ps}) for the product of two Chebyshev
functions of the second kind allows us as an immediate consequence
to find differential equations satisfied by Chebyshev exponents.
\begin{theorem}
If $y = x^{\cp k}$, then
\begin{equation}
(x^{2}-4)(y')^{2}-k^{2}(y^{2}-4)=0 \label{eq:ford}
\end{equation}
\end{theorem}
Proof:
From (\ref{eq:evefact}) we obtain
\[(x^{2}-4)\Sc_{k}^{2}(x) = x^{\cp 2k}-2=(x^{\cp k})^{2}-4.\] 
Substituting $y'/k$ for $\Sc_{k}(x)$, we obtain the theorem.
$\Box$

From this nonlinear equation of order one we can now find a linear
equation of order two.
\begin{theorem}
If $y = x^{\cp k}$, then
\begin{equation}
(x^{2}-4)y''+xy'-k^{2}y = 0 \label{eq:sord}
\end{equation}
\end{theorem}
Proof:
Take the derivative of (\ref{eq:ford}) and divide by $2y'$.
$\Box$

These equations are familiar as the differential equations satisfied
by Chebyshev polynomials for integral $k$; here we are assuming
$k$ may be any complex number.

We now may find a recursive expression for higher derivatives.
\begin{theorem}
If $y = x^{\cp k}$, then
\begin{equation}
(x^{2}-4)y^{(n+2)}+(2n+1)xy^{(n+1)}+(n^{2}-k^{2})y^{(n)} \label{eq:sigh}
\end{equation}
\end{theorem}
By an induction hypothesis, we may assume
\[(x^{2}-4)y^{(n+1)}+(2n-1)xy^{(n)}+((n-1)^{2}-k^{2})y^{(n-1)}=0.\]
Taking the derivative, we obtain the theorem.
$\Box$

Since (\ref{eq:sord}) is a second-order linear equation, it must
have a second linearly independent solution. When $k$ is not an
integer we have the following.
\begin{proposition}
\label{2sol}
For any $k$, $y = (-x)^{\cp k}$ satisfies the differential
equation (\ref{eq:sord}). When $k$ is not an integer, it is
linearly independent of $x^{\cp k}$.
\end{proposition}
Proof:
Substituting $-x$ for $x$ and $-y'$ for $y'$ in (\ref{eq:sord})
leaves it unchanged, and hence $(-x)^{\cp k}$ is a solution.

By Taylor's theorem, we may express $x^{\cp k}$ in a power series
of radius of convergence two around $x=0$, since the nearest 
singularity is at $x=-2$. We then have
\[x^{\cp k} = 0^{\cp k}{\cal E}(x) + k \Sc_{k}(0) {\cal O}(x),\]
where ${\cal E}$ is an even function, and ${\cal O}$ is an odd
function. If $k$ is an odd integer, then $0^{\cp k}=0$, whereas
if it is an even integer, $k \Sc_{k}(0)=0$; otherwise it is
non-zero. Hence assuming $k$ is not an integer,
\[(-x)^{\cp k} = 0^{\cp k}{\cal E}(x) - k \Sc_{k}(0) {\cal O}(x)\]
is linearly independent of $x^{\cp k}$.
$\Box$

This theorem is illuminating in connection with the problem
of defining the branches of the Chebyshev exponent, 
since any solution of the equation (\ref{eq:sord}) for non-integral
$k$ must be expressible at least locally
as a linear combination of $x^{\cp k}$ and $(-x)^{\cp k}$.

We may approach the matter in another way, which will also cover
the case where $k$ is a non-zero integer.
\begin{proposition}
Let 
\[y= \sqrt{x^{2}-4}\Sc_{k}(x).\]
Then $y$ satisfies the second order differential equation
(\ref{eq:sord}) for the Chebyshev exponent, namely
\[(x^{2}-4)y''+xy'-k^{2}y = 0.\]
\end{proposition}
Proof:
Let $w=\Sc_{k}(x)$. From (\ref{eq:sigh}) we have
\[(x^{2}-4)w''+ 3xw' + (1-k^{2})w = 0.\]
By differentiating, we find that
$$\sqrt{x^{2}-4}w = y,$$
$$\sqrt{x^{2}-4}w' = y'-\frac{xy}{x^{2}-4},$$
$$\sqrt{x^{2}-4}w'' = y'' - \frac{y+2xy'}{x^{2}-4} + 
\frac{3x^{2}y}{(x^{2}-4)^{2}}.$$
Substituting, we obtain (\ref{eq:sord}). 
$\Box$

\begin{proposition}
For all nonzero $k$, $\sqrt{x^{2}-4}\Sc_{k}(x)$ is linearly
independent of $x^{\cp k}$
\end{proposition}
Proof:
The Chebyshev power function $x^{\cp k}$ is ramified only at $x=-2$,
whereas $\sqrt{x^{2}-4}\Sc_{k}(x)$ ramifies when $x=2$.
$\Box$

Note that this second solution can be expressed in other ways,
up to a possible change of sign, for we have 
\[\sqrt{x^{2}-4}\Sc_{k}(x)= \pm \sqrt{x-2}\Uu_{2k}(x)= \pm \sqrt{x^{\cp 
2k}-2};\] 
for most purposes $\sqrt{x-2}\Uu_{2k}(x)$ is the preferred version.

It should also be noted that we could use instead the function
$\sqrt{x^{2}-4}\Sc_{k}(x)/k$ in all cases, for we have
the following limit.
\begin{proposition}
We have
\[\lim_{k \rightarrow 0}\frac{1}{k}\sqrt{x^{2}-4}\Sc_{k}(x)
= \pm 2 \log_{\cp } x.\]
\end{proposition}
Proof:
Let $x = z+z^{-1}$.
From the definition of $\Sc_{k}(x)$, we then have
\[\sqrt{x^{2}-4}\Sc_{k}(x) = \pm (z^{k}-z^{-k}).\]
 By l'H\^{o}pital's rule,
\[\lim_{k \rightarrow 0}\frac{z^{k}-z^{-k}}{k} = 2 \log z = \pm 2 
\log_{\cp} x,\]
and hence the proposition.
$\Box$ 

\subsection{Power Series Expansions}

Despite any appearances to the contrary, 
$x^{\cp k}$ is analytic in a neighborhood of $x=2$.
In fact, in analogy with the binomial
expansion of $(1+x)^{k}$, we have the following.
\begin{theorem}
For any complex $k$ and any $x$ such that $|x| < 4$, 
$(2+x)^{\cp k}$ can be expressed as a power series around $x=0$ by
\[ (2+x)^{\cp k} = 
2 + k \sum_{i=1}^{\infty}{k+i-1 \choose 2i-1}{x^i \over i} \]
\begin{equation}
\label{pows2}
=2+k^{2}x+\frac{2}{4!}k^{2}(k^{2}-1)x^{2}+
\frac{2}{6!}k^{2}(k^{2}-1)(k^{2}-4)x^{3}+ \cdots .
\end{equation}
\end{theorem}
Proof:
From ({\ref{eq:sigh}), we have that when $x=2$, if $y=x^{\cp x}$,
then
\[y^{(n+1)}=\frac{k^{2}-n^{2}}{4n+2}y^{(n)}.\]
From this we deduce that if
\[y=\sum_{n=0}^{\infty}a_{n}x^{n}\]
then 
\[a_{n+1} = \frac{1}{n+1}\frac{k^{2}-n^{2}}{4n+2}a_{n}
=\frac{k^{2}-n^{2}}{(2n+1)(2n+2)}a_{n}.\]
Since $a_{0}=2$, the series follows by induction.

Since 
\[\lim_{n \rightarrow \infty} \frac{a_{n+1}}{a_{n}} = -\frac{1}{4},\]
the series converges absolutely when $|x|<4$, as expected,
and does not converge anywhere on the circle $|x|=4$.
$\Box$

We now also have
\begin{corollary}
\label{S-expand}
For $|x| < 4$, we have
\begin{equation}
\Sc_{k}(2+x) = \sum_{i=0}^{\infty}{k+i \choose 2i+1}x^{i}
\end{equation}
\end{corollary}
Proof: Termwise differentiation.
$\Box$

The reader familiar with hypergeometric functions will not
be surprised to learn that
the above power series expansions are special cases of the
classical hypergeometric series of Gauss.
\begin{theorem}
Let 
\[F(a,b,c,x) = \sum_{i=0}^{\infty}\frac{{a+i-1 \choose i}
{b+i-1 \choose i}}{{c+i-1 \choose i}}x^{i}\]
be the hypergeometric series of Gauss, convergent (at least) when
$|x|<1$. Then for $|x| < 4$,
\begin{equation}
(2+x)^{\cp k}=2F(k,-k, \frac{1}{2},-\frac{x}{4})
\end{equation}
\end{theorem}
Proof:
Comparison of the two series shows they are identical.
$\Box$

With a corresponding choice of branch, we might now define
the Chebyshev exponent as
\begin{equation}
x^{\cp k} = 2 F(k,-k,\frac{1}{2},\frac{2-x}{4}) \label{eq:hyper}
\end{equation}
where $F$ is now the classical hypergeometric function, 
defined even where the hypergeometric series does not converge by 
analytic continuation.  We might then use the varied and sometimes 
complicated properties of the hypergeometric function to produce
a variety of formulas and series expansions. We will not pursue
the matter, as only a few series will be important for later sections.

We will also note that one property of the hypergeometric function is 
that its derivatives to all orders are also hypergeometric functions.  
This gives us a formula for the higher derivatives of Chebyshev 
powers.
\begin{proposition}
For any positive integer $n$, we have
\begin{equation}
(x^{\cp k})^{(n)} = (n-1)!k {k+n-1 \choose 2n-1}
F(n+k,n-k,n+\frac{1}{2},\frac{2-x}{4})
\end{equation}
\end{proposition}
Proof: 
It is immediate from the power series expansion that
\[F'(a,b,c,x) = \frac{ab}{c}F(a+1,b+1,c+1,x).\]
Applying this to (\ref{eq:hyper}) gives us the theorem.
$\Box$

We also have a nice expression for $\Uu_{k}$ in terms of
a power series expansion around $x=2$.
\begin{theorem}
For any complex $k$ and for $|x| < 4$, we have
\begin{equation}
\Uu_{k}(2+x)  = k \sum_{i=0}^{\infty}{(k-1)/2+i \choose 
2i}\frac{x^{i}}{2i+1} \label{eq:u-series}
\end{equation}
\begin{equation}
= k + \frac{k(k^{2}-1^{2})}{3!}(\frac{x}{4})+
\frac{k(k^{2}-1^{2})(k^{2}-3^{2})}{5!}(\frac{x}{4})^{2}+\cdots
\end{equation}
\end{theorem}
Proof:
From \ref{S-expand}, we have
\[\Uu_{k}(2+x) = \sum_{i=0}^{\infty}({(k-1)/2+i \choose 2i+1}+{(k+1)/2+i \choose 
2i+1})x^{i}.\]
For any non-negative integer $n$, we define a polynomial $p_{n}$ of 
degree $2n+1$ by setting 
\[p_{n}(x) = {(x-1)/2+n \choose 2n+1}+{(x+1)/2+n \choose 2n+1}.\]
We then have $p_{n}(m) = 0$ for any odd integer $-2n < m < 2n$, 
and we also have $p_{n}(0) = 0$. 
Hence we have $2n+1$ distinct roots for a
polynomial of degree $2n+1$, which must therefore all be simple roots.
We also have a coefficient on the leading term of 
$1/(4^{n}(2n+1)!)$.
Since these conditions determine $p_{n}$, we have
\[p_{n}(x) = \frac{x}{4^{n}(2n+1)}{(x-1)/2+n \choose 2n},\]
since they are also plainly true of the right-hand side.
Hence the $i$-th term of the series is $p_{i}(k)x^{i}$, and
we have the theorem.
$\Box$

\begin{corollary}
For an appropriate choice of branch cut, we have
\[\Uu_{k}(x) = 
k F(\frac{1+k}{2},\frac{1-k}{2},\frac{3}{2},\frac{2-x}{4}).\]
\end{corollary}
Proof:
Comparison of the power series expansion around $x=2$.
$\Box$

We may use the expansion for $\Uu$ in expressing $x^{\cp k}$ in power 
series; we began with the expansion around $x=0$.
\begin{theorem}
\label{pow-zero}
The power series for $x^{\cp k}$ expanded around $x=0$ and
valid for $|x| < 2$ can be obtained by expanding
\begin{equation}
\label{eq:pow-zero}
x^{\cp k}=\cos (\frac{\pi}{2}k) (2-x^{2})^{\cp k/2} + 
\sin (\frac{\pi}{2}k) x \Uu_{k}(2-x^{2})
\end{equation}
\end{theorem}
Proof:
Substituting $x=0$ into (\ref{eq:sigh}), we obtain
\[y^{(n+2)} = \frac{n^{2}-k^{2}}{4}y^{(n)},\]
so that if $a_{n}$ is the coefficient of $x^{n}$, we have
\[a_{n+2} =\frac{n^{2}-k^{2}}{4(n+1)(n+2)}a_{n}.\]
Using this, we find the even part of the power series expansion is
\[0^{\cp k}(1-\frac{k^{2}}{2!}(\frac{x}{2})^{2}
+\frac{k^{2}(k^{2}-2^{2})}{4!}(\frac{x}{2})^{4}
-\frac{k^{2}(k^{2}-2^{2})(k^{2}-4^{2})}{6!}(\frac{x}{2})^{6}+ \cdots ).\]
Noting that
\[\frac{0^{\cp k}}{2} = \frac{i^{k}+i^{-k}}{2} = \cos \frac{\pi}{2}k,\]
and comparing power series expansions, we find that this is
\[\cos \frac{\pi}{2}k (2-x^{2})^{\cp k/2}.\]

In the same way, we find that the odd part of the power series
expansion is
\[k \Sc_{k}(0)x(1 - \frac{k^{2}-1^{2}}{3!}(\frac{x}{2})^{2}
+\frac{(k^{2}-1^{2})(k^{2}-3^{2})}{5!}(\frac{x}{2})^{4} + \cdots,\]
by noting that
\[\Sc_{k}(0) = \frac{i^{k}-i^{-k}}{i-i^{-1}} = \sin \frac{\pi}{2}k,\]
and comparing power series expansions with
\[\sin (\frac{\pi}{2}k) x \Uu_{k}(2-x^{2})\]
we see that they are identical.
$\Box$

Of course, having established this identity when $|x| \le 2$ we
may continue it analytically. We then also get the following 
identities, valid everywhere.
\begin{corollary}
When $k$ is not an odd integer, we have
\[(2-x^{2})^{\cp k/2} = \frac{x^{\cp k} + (-x)^{\cp k}}{2 
\cos(\pi k/2)},\]
while if $k$ is not an even integer we have
\[x \Uu_{k}(2-x^{2})=\frac{x^{\cp k}-(-x)^{\cp k}}{2 
\sin(\pi k/2)}.\]
\end{corollary}

We may also use this identity to derive the
Puiseux series expansion (in powers of $\sqrt{x+2} = \sqrt[\cp]{x}$)
around $x=-2$.
\begin{theorem}
The Puiseux series expansion of $x^{\cp k}$ around the ramified
place $x=-2$ is given by expanding each part of
\[x^{\cp k} = \cos (\pi k) (-x)^{\cp k}
+ \sin (\pi k) \sqrt{x+2} \Uu_{2k}(-x)\]
in power series.
\end{theorem}
Proof:
We may write (\ref{eq:pow-zero}) as
\[x^{\cp k}=\cos (\frac{\pi}{2})k (-x^{\cp 2})^{\cp k/2} + 
\sin (\frac{\pi}{2}k) x \Uu_{k}(-x^{\cp 2}).\]
If rewrite $x^{\cp k}$ in this expression as $(\sqrt[\cp]{x})^{\cp 2k}$,
we obtain
\[x^{\cp k}=\cos (\pi k) (-x)^{\cp k} + 
\sin (\pi k) \sqrt[\cp]{x} \Uu_{2k}(-x).\]
$\Box$

\section{Polynomial Equations}

\subsection{The $n$ Chebyshev $n$-th Roots}

In any field containing $n$-th roots of unity, if we know an element 
$r_{0}$ such that $r_{0}^{n}=t$, it is easy to find the other $n-1$ 
roots of $x^{n}-t=0$; for we simply take a primitive $n$-th root of 
unity $\zeta$, and the other roots will be $r_{i}=\zeta^{i}r_{0}$.  
The situation is similar though a bit more complicated for Chebyshev 
roots.

\begin{definition}
Let $n$ be a positive integer, and let $r$ be $n$ elements of
an algebraically closed field of characteristic not dividing $n$,
indexed by residue classes modulo $n$. Let $\zeta$ be a primitive
$n$-th root of unity. Noting if $i$ is a residue class modulo
$n$ that $\zeta^{n}$ is well-defined, we will call $r$ an
{\em indexed set of $n$-th Chebyshev roots of $t$} if 
$r_{0}=u+u^{-1}$, $r_{0}^{\cp n} =t$, and
\[r_{i} = \zeta^{i}u + \zeta^{-i}u^{-1}.\]
\end{definition}

\begin{proposition}
Let $r$ be an indexed set of $n$-th Chebyshev roots of $t$.
Then the $n$ roots of $x^{\cp n}-t=0$ are in fact given by $r$.
\end{proposition}
Proof:
\[r_{i}^{\cp n} = u^{n}+u^{-n}=r_{0}^{\cp n}=t.\]
$\Box$

\begin{proposition}
Let $r$ be an indexed set of $n$-th Chebyshev roots of $t$.
Then 
\[r_{k+i}+r_{k-i} = \mu^{\cp i}r_{k},\]
where $\mu = \zeta+\zeta^{-1}$ is a primitive $n$-th Chebyshev root of
two corresponding to the primitive $n$-th root of unity $\zeta$.
\end{proposition}
Proof:
By the previous theorem,
\[r_{k+i}+r_{k-i} = 
(\zeta^{k+i}+\zeta^{k-i})u + (\zeta^{-k-i}+\zeta^{i-k})u^{-1}\]
\[= (\zeta^{i}+\zeta^{-i})\zeta^{k}u + 
(\zeta^{i}+\zeta^{-i})\zeta^{-k}u^{-1}
= \mu^{\cp i}r_{k}.\]
$\Box$

We may express all of the roots of in terms of just one of them, by
means of a quadratic equation.
\begin{proposition}
Let $r$ be an indexed set of $n$-th Chebyshev roots of $t$, with
$\mu$ as before.
Then we have that $r_{k+i}$ and $r_{k-i}$ are the two roots of
\[x^{2}-\mu^{\cp i}r_{k}x + r_{k}^{\cp 2}+\mu^{\cp 2i} = 0.\]
\end{proposition}
Proof:
The trace term (the coefficient of $x$) is given by the previous
proposition. The constant term we may find by multiplying out
\[r_{k+i}r_{k-i} = (\zeta^{k+i}u + \zeta^{-k-i}u^{-1})
(\zeta^{k-i}u + \zeta^{i-k}u^{-1}) \]
\[=\zeta^{2k}u^{2} + z^{-2k}u^{-2} + \zeta^{2i}+\zeta^{-2i}
=r_{k}^{\cp 2} + \mu^{\cp 2i}.\]
$\Box$

We may also express all of the roots as a linear combination of
two of them.
\begin{proposition}
Let $\mu$ be as before, 
let $r$ be an indexed set of $n$-th Chebyshev roots of $t$, and
let $r_{i}$ and $r_{j}$ be two of the roots. Then if $e=j-i$ and if
$k$ is an integer, we have that
\begin{equation}
r_{i+ke}=\Sc_{k}(\mu^{\cp e})r_{j}-\Sc_{k-1}(\mu^{\cp e})r_{i}.
\end{equation}
\end{proposition}
Proof:
The proposition is clearly true when $k=0$ or $k=1$, and by induction
we may assume it to be true up to $k-1$. Then
\[r_{i+ke}+r_{i+(k-2)e} = \mu^{\cp e}r_{i+(k-1)e}.\]
Hence
\[r_{i+ke}= \mu^{\cp e}(\Sc_{k-1}(\mu^{\cp e})r_{j}
-\Sc_{k-2}(\mu^{\cp e})r_{i})-\Sc_{k-2}(\mu^{\cp e})r_{j}
+\Sc_{k-3}(\mu^{\cp e})r_{i}\]
\[=(\mu^{\cp e}\Sc_{k-1}(\mu^{\cp e})-S_{k-2}(\mu^{\cp e})r_{j}
-(\mu^{\cp e}\Sc_{k-2}(\mu^{\cp e})-\Sc_{k-3}(\mu^{\cp e})r_{i}\]
\[=\Sc_{k}(\mu^{\cp e})r_{j}-\Sc_{k-1}(\mu^{\cp e})r_{j}.\]
The theorem for negative values of $k$ now follows immediately.
$\Box$

\begin{theorem}
\label{inset}
Let $r$ be an indexed set of $n$-th Chebyshev roots of $t$, and let 
$i$ and $j$ be indices such that $e=j-i$ is prime to $n$.  This last 
condition means that any element of the residue class $e$ is prime to 
$n$, and entails $e \ne 0$, so that $i \ne j$.  If $\mu$
is as before, then if the base field 
is of characteristic 0, any root $r_{l}$ of $t$ can be written as a 
$\Z[\mu]$-linear combination of $r_{i}$ and $r_{j}$; while if it is of 
characteristic $p$, with $n$ not divisible by $p$, then it is an 
$F_{p}(\mu)$-linear combination of $r_{i}$ and $r_{j}$, where $F_{p}$ 
is the field of $p$ elements.
\end{theorem}
Proof:
That $e$ is prime to $n$ means that $e$ is invertible, which is to
say, $ke$ for values of $k$ from 0 to $n-1$ traverse a complete set
of residue classes modulo $n$. Hence $i+ke$ traverses a complete set
of residue classes, and so $r_{i+ke}$ is a complete set of Chebyshev
roots of $t$, taken in another order. 

In characteristic 0,
$\Q(\mu)$ is
cyclic and since $e$ is prime to $n$, $\Q(\mu^{\cp e}) = \Q(\mu)$,
and any polynomial in $\mu^{\cp e}$ can also be expressed as
a polynomial (of degree less than $\phi(n)/2$) in $\mu$. Hence any
root is a $\Z[\mu]$-linear combination of $r_{i}$ and $r_{j}$. The
theorem in characteristic $p$ follows immediately on reduction
modulo $p$.
$\Box$

\subsection{Branches of the Chebyshev Radical}

When considering ordinary $n$-th roots, it is of course of
considerable utility to have in mind the representation of
complex numbers in terms of polar coordinates. If we wish for 
something similar in the case of Chebyshev $n$-th roots, 
we may began by considering the curves which arise by setting
the real or imaginary part of $\log_{\cp} z$ to a constant
value.

If $\log_{\cp} z = r + i \theta$, then
\[z = \exp_{\cp}(r+i \theta) =
2 \cosh r \cos \theta +2 \sinh r \sin \theta.\]
If $x = 2 \cosh r \cos \theta$ and $y = 2 \sinh r \sin \theta$, then
\begin{eqnarray}
\frac{x^{2}}{4 \cosh^{2} r} + \frac{y^{2}}{4 \sinh^{2} r} = 1 
\label{eq:ell} \\
\frac{x^{2}}{4 \cos^{2} \theta} - \frac{y^{2}}{4 \sin^{2} \theta} = 1 
\label{eq:hyp}
\end{eqnarray}
as we may check by substituting.

Hence if $r$ is a constant, (\ref{eq:ell}) gives the equation of an 
ellipse which we may take as the analog of a circle centered on the 
origin in polar coordinates.  The analog of a ray from the origin of 
constant angular argument $\theta$ would be either the top or the 
bottom part of one branch of the hyperbola (\ref{eq:hyp}).  To get the 
analog of a line through the origin, we take a complete branch of the 
hyperbola, and hence half a hyperbola instead of one quarter of one.

We may view the matter more schematically by representing any complex 
number $z$ which is not a real number of absolute value less than or 
equal to two by the polar coordinates $r$ and $\theta$, where 
$\log_{\cp} z = r + i \theta$.  Real $z$ between $-2$ and $2$ 
correspond to the origin; however, since $\theta$ varies, we may 
``blow up'' the origin into a circle which doubly covers the line 
segment $-2 \le z \le 2$ by taking $\theta$ in the range $-\pi < \theta 
\le \pi$, so that we have a complete circle, and then identifying the 
point represented by $\theta$ with that represented by $-\theta$.

Bearing the above representations in mind, we may make the
following definition.
\begin{definition}
Let $t$ be any complex number, $n$ a positive integer, $l$
any integer, and
$z_{0} = \sqrt[\cp n]{t}$ the $n$-th Chebyshev root of $t$. If
\[\log_{\cp}z_{0} = r + i \theta_{0},\]
we define an indexed set of complex branches of the $n$-th
Chebyshev root by
\[z_{l}=\sqrt[\cp n<l>]{t} = \exp_{\cp}(r + i \theta_{l}),\]
where either $r \ge 0$ and
\begin{equation}
\frac{l \pi}{n} \le \theta_{l} < \frac{(l+1) \pi}{n} \label{eq:evecon}
\end{equation}
or $r > 0$ and
\begin{equation}
-\frac{(l+1) \pi}{n} \le \theta_{l} < -\frac{l \pi}{n} 
\label{eq:oddcon}
\end{equation}
obtains;
and such that
\[\theta_{l} \equiv \theta_{0} \pmod{\frac{2 \pi}{n}}.\]
\end{definition}

\begin{definition}
For any rational exponent $p/q$, we define $x^{\cp p/q<l>}$ as
\[x^{\cp p/q<l>} = (\sqrt[\cp q<l>]{x})^{p}.\]
\end{definition}

\begin{proposition} 
The functions $\sqrt[\cp n<i>]{}$ and $\sqrt[\cp n<j>]{}$
are identical if and only if
$i \equiv j \pmod {2n}$ or $i+j+1 \equiv 0 \pmod {2n}$.
\end{proposition}
Proof:
If $i \equiv j \pmod {2n}$, then (\ref{eq:evecon}) and 
(\ref{eq:oddcon}) differ by a multiple of $2 \pi$, which is
the period of $\exp_{\cp}$. Substituting $-l-1$ for $l$
in (\ref{eq:evecon}) gives us (\ref{eq:oddcon}), and vice-versa.
Putting these together, we get the two conditions of the theorem.
No more conditions are possible, since $i+j+1 \equiv 0 \pmod {2n}$
identifies pairs of congruence classes modulo $2n$, and hence
leads to $n$ classes, which may not be further identified since
we have $n$ branches.
$\Box$

\begin{theorem}
The $n$ functions defined above do in fact constitute $n$ distinct
branches of the $n$-th Chebyshev root.
\end{theorem}
Proof:
It is immediate that 
\[(\sqrt[\cp n<l>]{t})^{\cp n}=\exp_{\cp}(nr+i n\theta_{l}) =
\exp(nr+i n \theta_{0}+2 \pi m) = z_{0}^{\cp n} = t,\]
since $m$ is an integer; and hence each $\sqrt[\cp n<l>]{t}$ is an
$n$-th Chebyshev root of $t$.

Since $\exp_{\cp}$ is an even function, we are at liberty to replace 
$\log_{cp}$ by $-\log_{\cp}$, and so 
(\ref{eq:oddcon}) by the condition $r < 0$ and
\[\frac{l \pi}{n} < \theta_{l} \le \frac{(l+1) \pi}{n}.\]
We therefore have that an equivalent condition is merely
(\ref{eq:evecon}) with $r$ allowed to take any real value.

Multiplication by $n$ now gives a continuous bijection
between the pathwise connected region so defined
and the region (for all $r$)
\[l \pi < \theta \le (l+1) \pi, \]
and Chebyshev exponentiation by $\exp_{\cp}$ gives a continuous
bijection from this to the whole complex plane. Hence 
$x = \sqrt[n<l>]{t}$
gives a bijection from the complex plane to a pathwise connected
region of the complex plane, whose inverse map is the $n$-th Chebyshev
power map $t = x^{\cp n}$. 
Since there are $n$ non-overlapping
regions defined by these branches which together cover the complex
plane, we have defined a complete set of branches.
$\Box$

Since $r_{i} = \sqrt[\cp n<2i>]{t}$ gives us an indexed set of $n$-th 
Chebyshev roots of $t$, we could use Theorem \ref{inset} to write all 
of the roots in terms of two of them.  However, it seems better to 
refer back instead to the linear differential equation (\ref{eq:sord}) 
satisfied by Chebyshev powers, and to note that any solution must, at 
least locally, be a linear combination of $\sqrt[\cp n]{t}$ and 
$\sqrt[\cp n]{-t}$; we therefore seek to write the branches as a 
combination of these functions, and hence in terms of the principal 
branch.
\begin{theorem}
Let $\mu = \exp_{\cp}(\pi i / n)=\sqrt[\cp n]{-2}$. Then
for even integers $0 \le i < n$ we have
\begin{equation}
\sqrt[\cp n<i>]{t} = \Sc_{i+1}(\mu) \sqrt[\cp n]{t}
- \Sc_{i}(\mu) \sqrt[\cp n]{-t},
\end{equation}
while for odd integers $0 \le i < n$ we have
\begin{equation}
\sqrt[\cp n<i>]{t} = -\Sc_{i}(\mu) \sqrt[\cp n]{t}
+ \Sc_{i+1}(\mu) \sqrt[\cp n]{t}.
\end{equation}
\end{theorem}
Proof:
Let $F(t) = \sqrt[\cp n<i>]{t}$, and let $E(t)$ and $O(t)$ be the
even and odd parts of $\sqrt[\cp n]{t}$, so that
\[E(t) = \frac{\sqrt[\cp n]{t}+\sqrt[\cp n]{-t}}{2},\]
\[O(t) = \frac{\sqrt[\cp n]{t}-\sqrt[\cp n]{-t}}{2}.\]
Since $F$ satisfies (\ref{eq:sord}) and $E$ and $O$ form
a basis for the solutions, we have, at least in a neighborhood
of zero,
\[F(t) = \frac{F(0)}{E(0)}E(t) + \frac{F'(0)}{O'(0)}O(t).\]

We now find 
\begin{eqnarray*}
E(0) = \sqrt[\cp n]{0} = 2 \cos (\frac{\pi}{2n}), \\
F(0) = 0^{\cp (2i+1)/n} = 2 \cos (\frac{\pi}{2n}(2i+1)), \\
O'(0) = \frac{1}{n}\Sc_{1/n}(0) = \sin(\frac{\pi}{2n})/{n}, \\
F'(0) = \frac{1}{n \Sc_{n}(0^{\cp (2i+1)/n})} = 
(-1)^{i} \sin (\frac{\pi}{2n} (2i+1))/{n}.
\end{eqnarray*}
Hence we have
\begin{equation}
F(t) = (\cos (\frac{\pi}{2n}(2i+1))/\cos (\frac{\pi}{2n}))E(t)
+ (-1)^{i}(\sin (\frac{\pi}{2n}(2i+1)) /\sin (\frac{\pi}{2n}))O(t)
\end{equation}

We now write everything with a common denominator of
\[2 \sin(\frac{\pi}{2n}) \cos (\frac{\pi}{2n}) = \sin 
(\frac{\pi}{n}),\]
and apply the angle-sum trigonometrical identities
to write this in terms of $\sqrt[\cp n]{t}$ and $\sqrt[\cp n]{-t}$,
obtaining for even values of $i$
\begin{equation}
F(t) = (\sin (\frac{\pi}{n}(i+1))/\sin (\frac{\pi}{n}))
\sqrt[\cp n]{t} - (\sin (\frac{\pi}{n}i) / \sin 
(\frac{\pi}{n})) \sqrt[\cp n]{-t},
\end{equation}
and for odd values of $i$
\begin{equation}
F(t) = -(\sin (\frac{\pi}{n}i) / \sin (\frac{\pi}{n}))
\sqrt[\cp n]{t} + (\sin (\frac{\pi}{n}(i+1)) / \sin 
(\frac{\pi}{n})) \sqrt[\cp n]{-t}.
\end{equation}
Putting this into the more algebraic language we used for general
fields, we obtain the identities in the theorem.

All that remains is to check where these identities are valid;
but since they mesh together in the correct way along the branch
cuts from $(2, \infty)$ and $(-\infty, -2)$ they are in
fact valid everywhere.
$\Box$

\subsection{Polynomial Equations and Radicals}

Perhaps the best known of all the uses to which the extraction of 
roots has been put is in the solution of algebraic equations in terms 
of radicals.  This has been productive not so much for what it 
accomplished, but for what it failed to accomplish.  It failed to 
express real roots in terms of real numbers, and in so doing it lead 
to the discovery of complex numbers.  It failed to express the roots 
of all polynomials with integer coefficients, and in so doing it led 
to the discovery of Galois theory.  With the notable exception of the 
quadratic formula, it has also failed to be of much practical use.

Some of these failures can be ameliorated by the simple expedient of 
employing Chebyshev radicals in the place of ordinary ones.  In 
particular, while expressions for roots in terms of Chebyshev radicals 
will still grow more complicated with increasing degree, they are 
nearly always less complicated.  It is in fact rather remarkable with 
what persistence people have insisted on radicals as the only 
canonical family of algebraic functions, given the problems that this 
sometimes entails.

In practice, a solution in radicals usually means an expression for
the roots of a polynomial in terms of the function $\sqrt[n]{z}$.
In a general algebraic context, it is taken to mean a solution in
terms of a tower of extensions, each of which is given by an
$n$-th root, for various $n$. The following definition is standard.
\begin{definition}
Let $P$ be a univariate polynomial with coefficients in a field
$F_{0}$, and
suppose that $P$ may be completely factored into linear
factors in a field $F_{n}$, where we have for each $i$ with
$0 < i \le n$ that
$F_{i} = F_{i-1}(r_{i})$ with $r_{i}^{m_{i}} \in F_{i-1}$. We then
say that $P$ is {\em solvable in radicals} relative to 
the field $F_{0}$.
\end{definition}

It is of course obvious what the Chebyshev version of this should
be.
\begin{definition}
Let $P$ and $F_{0}$ be as before, and 
suppose that $P$ may be completely factored into linear
factors in a field $F_{n}$, where we have for each $i$ with
$0 < i \le n$ that
$F_{i} = F_{i-1}(r_{i})$ with $r_{i}^{\cp m_{i}} \in F_{i-1}$. We then
say that $P$ is {\em solvable in Chebyshev radicals} relative
to the field $F_{0}$.
\end{definition}

We now may prove the following fundamental theorem.
\begin{theorem}
A polynomial with coefficients in a field whose characteristic is not 
two is solvable in radicals if and only if it is solvable in Chebyshev 
radicals.
\end{theorem}
Proof:
Let us suppose that a polynomial $P$ has been solved in radicals
with respect to a field $F_{0}$. We wish to see if it can also
therefore be solved in Chebyshev radicals. Since the
field $F_{n}$ is obtained as a tower of radical extensions over
$F_{0}$, it suffices to consider whether any polynomial of the
form $x^{n}-t$, with $t$ in a field $K$ of characteristic equal to
that of $F_{0}$, may be solved in terms of Chebyshev radicals.
Because of the composition property of exponents, we may 
reduce to the case where the exponent $x^{q}-t$ is a prime $q$;
and since in characteristic $p$, $x^{p}=x^{\cp p}$, we may further 
assume without loss of generality
that $q$ is not equal to the characteristic.

We wish therefore to create a tower of Chebyshev radical extensions 
which will factor the polynomial $x^{q}-t$, where $t \in K_{0}$, and 
where we may assume $t \ne 0$ since we have trivially $0 \in K_{0}$.

We first set $K_{1} = K_{0}(\mu)$, where $\mu \ne 2$ and
$\mu^{\cp q} = 2$. Then we have $K_{2}=K_{1}(\lambda)$, where
$\lambda^{\cp 2}=\mu^{\cp 2}-4$. Now set $\zeta = (\mu + \lambda)/2$.
The pair of equations $\mu^{2}+\lambda^{2}=4$, $\mu+\lambda=2$
have solutions $\mu=2, \lambda = 0$ and $\mu=0, \lambda = 2$.
The first solution is ruled out since we have assumed
$\mu \ne 2$, and the second is ruled out since $0^{\cp q} \ne 2$
for any prime $q$. Hence $\zeta \ne 1$. However, we have
\[ (\frac{\mu + \lambda}{2}) (\frac{\mu - \lambda}{2})
= \frac{\mu^{\cp 2}-\lambda^{\cp 2}}{4} = 1,\]
and so $\zeta^{-1} = (\mu-\lambda)/2$. Also, 
\[\zeta^{q} + \zeta^{-q} = (\zeta+\zeta^{-1})^{\cp q}
=\mu^{\cp q} = 2,\]
and so 
\[\zeta^{2q}-2\zeta^{q}+1 = (\zeta^{q}-1)^{2}=0.\]
It follows that $\zeta$ is a $q$-th root of unity not equal to one.

We now set $K_{3} = K_{2}(s)$, where $s^{\cp q}=t+t^{-1}$.  Finally, 
we set $K_{4}=K_{3}(r)$, where $r^{\cp 2} = s^{\cp 2}-4$.

We now proceed just as before, noting first that
\[(\frac{s+r}{2})(\frac{s-r}{2}) = \frac{s^{\cp 2}-r^{\cp 2}}{4} = 1,\]
and so if $w=(s+r)/2$, $w^{-1} = (s-r)/2$. Then we have
\[w^{q}+w^{-q}=(w+w^{-1})^{\cp q}=s^{\cp q} = t+t^{-1},\]
and so
\[w^{2q}-(t+t^{-1)})w^{q}+1= (w^{q}-t)(w^{q}-t^{-1}) = 0.\]
Hence either $w^{q}=t$ or $(w^{-1})^{q}=t$. We now may factor
$x^{q}-t$ by taking its roots to be $\zeta^{i}w$ or $\zeta^{i}w^{-1}$,
depending on which of $w$ or $w^{-1}$ has $q$-th powers equal to $t$.

To show the contrary direction, we need to create a tower of
ordinary radical extensions which will factor $x^{\cp q}-t$, since
we may make precisely the same reductions as before, and assume
without loss of generality that $q$ is a prime not equal to the
characteristic of $F_{0}$.

We began by setting $K_{1}= K_{0}(\zeta)$, where $\zeta^{q}=1$
but $\zeta \ne 1$. Now set $K_{2} = K_{1}(s)$, where
$s^{2} = t^{2}-4$. Finally, set $K_{3}=K_{2}(r)$, where
$r^{q}=(s+t)/2$. Noting as before that if $r^{q} = (s+t)/2$ we
have $r^{-q} = (s-t)/2$, we have
\[(\zeta^{i}r + \zeta^{-i}r^{-1})^{\cp q} = r^{q}+r^{-q} = s,\]
which therefore gives us the $q$ roots of $x^{\cp q}-t$.
$\Box$

In characteristic two, Chebyshev radicals turn out to solve
a wider class of polynomial equations than do ordinary radicals.
We first prove the following lemma.
\begin{lemma}
Let $F$ be a field of characteristic two, and let $a \ne 1$ be
an element of $F$. If $b^{\cp 3} = a/(a+1)^{3}$ and if
$b \ne a/(a+1)$, then if we set $c=(a+1)b$ then $c$ satisfies
the equation
\[c^{2}+ac+1 = 0.\]
\end{lemma}
Proof:
Since 
\[(\frac{c}{a+1})^{\cp 3} = (\frac{c}{a+1})^{3}+(\frac{c}{a+1})=
\frac{a}{(a+1)^{3}},\]
we have on multiplying through by $(a+1)^{3}$ that
\[c^{3}+(a+1)^{2}c + a = 0.\]
Since $b \ne a/(a+1)$ we have $c \ne a$; hence we may divide by
$c+a$ and obtain
\[c^{2}+ac+1 = 0.\]
$\Box$

We now have the following.
\begin{theorem}
If $P$ is a polynomial with coefficients in a field of characteristic
two, then if $P$ can be solved in radicals it can be solved
in Chebyshev radicals, but not conversely.
\end{theorem}
Proof:
As before, we may reduce the problem to factoring $x^{q}+t$, where $q$ 
is an odd prime, over a field obtained as a tower of Chebyshev radical 
extensions of $K_{0}$. We first show that we may introduce $q$-th
roots of unity.

If $\zeta \ne 1$ and $\zeta^{3} = 1$, then 
$\zeta^{\cp 5} = \zeta (\zeta^{2} + \zeta +1)^{2}$, and so
we may define $\zeta$ as a Chebyshev radical by setting
$\zeta^{\cp 5} = 0$, $\zeta \ne 0$. In case $q=3$, we define
$K_{1} = K_{0}(\zeta)$, and $K_{2} = K_{1}$.

Now suppose $q \ne 3$ and we want to express $\zeta$ such that
$\zeta^{q} = 1$ while $\zeta \ne 1$ by a Chebyshev radical 
extension. Let $\mu^{q}=0$, while $\mu \ne 0$. We now set
$F_{1} = F_{0}(\mu)$, and let $\lambda^{\cp 3} = \mu/(\mu+1)^{3}$
with $\lambda \ne \mu$.
Since $q \ne 3$, $\mu \ne 1$ and this is well-defined. We
now set $K_{2} = K_{1}(\lambda)$.
Invoking the lemma, if $\zeta = (\mu+1) \lambda$ then
$\zeta^{2}+\mu \zeta + 1 = 0$, and hence $\mu = \zeta + \zeta^{-1}$,
from which it follows that $\zeta^{q}=1$ while $\zeta \ne 1$,
and hence $K_{2}$ contains the $q$-th roots of unity.

We now set as before $K_{3} = K_{2}(s)$, where $s^{\cp q} = t + 
t^{-1}$. Next we define $r$ by $r^{\cp q} = s/(s+1)^{3}$,
$r \ne s$, and set $K_{4} = K_{3}(r)$. If $w = (s+1)r$, we invoke
the lemma and conclude that $w^{2}+sw + 1 = 0$. Hence just
as before we have
\[w^{q}+w^{-q}=(w+w^{-1})^{\cp q}=s^{\cp q} = t+t^{-1},\]
and so
\[w^{2q}-(t+t^{-1)})w^{q}+1= (w^{q}-t)(w^{q}-t^{-1}) = 0;\]
so that either $w^{q} =t$ or $(w^{-1})^{q} = t$.
Again we conclude that with the roots of unity in place, we may
completely factor $P$ over the field $K_{4}$

Let us now suppose that $F$ is a field of characteristic two 
containing the $q$-th roots of unity for every odd prime $q$, that
$t$ is an element transcendental over $F$, and let $K_{0} = F(t)$.  From 
the lemma, the polynomial $x^{2}+tx+1$ may be solved in Chebyshev 
radicals over $K_{0}$.  We will show it cannot be solved in ordinary 
radicals over $K_{0}$.

Suppose $K_{n}$ is a field obtained as a tower of radical extensions 
of prime degree over $K_{0}$ in which $x^{2}+tx+1$ may be factored.  If 
$K_{n}/K_{n-1}$ is a radical extension, obtained by adjoining a square 
root, then the polynomial can be factored over $K_{n-1}$, since 
$x^{2}+tx+1$ is separable.  On the other hand, if $K_{n} = K_{n-1}(r)$ 
where $r^{q} \in K_{n-1}$ and $q$ is an odd prime, then either $r \in 
K_{n-1}$ and so $K_{n} = K_{n-1}$, or $K_{n}$ is cyclic of degree $q$ 
over $K_{n-1}$, in which case $x^{2}+tx+1$ must already be 
factorable over $K_{n-1}$, since there exists no subextension.  
Iterating this process, we conclude that $x^{2}+tx+1$ is factorable 
over $K_{0}$, which is a contradiction.
$\Box$

\subsection{Galois Theory}

Before discussing how to go about solving polynomial equations in 
Chebyshev radicals, we will recall for the reader of some of the basic 
facts of Galois theory; the reader who wishes something more than this 
should consult any of the standard textbooks.

If $P$ is a polynomial with coefficients in $F$ for which we seek a 
solution in radicals, we may always first factor $P$ over $F$, and 
hence we may without loss of generality assume $P$ is irreducible over 
$F$.  In characteristic $p$, any irreducible $P$ may be written as 
$P(x) = Q(x^{q}) = Q(x^{\cp q})$, where $q = p^{n}$ is a power of $p$ 
and where $Q$ is irreducible and separable; that is, $(Q, Q') = 1$, so 
that $Q$ has no repeated factors.  Hence we may without loss of 
generality assume that $P$ is both irreducible and separable, since 
the $q$-th root extension is both an ordinary and a Chebyshev radical.

If $P$ is separable, its roots when adjoined to $F$ determine a {\em 
Galois extension} $K/F$.  The automorphisms of $K$ fixing $F$ 
constitute finite group called the {\em Galois group} of the extension 
$K/F$; this group acts faithfully on the roots of $P$, giving a 
faithful permutation representation of the Galois group as a 
permutation group, which is transitive if $P$ is irreducible.  

If the Galois group $G$ has a normal subgroup $H_{1}$ such that the 
quotient group $G/H_{1}$ is cyclic of prime degree, and $H_{1}$ in turn 
has a normal subgroup $H_{2}$ such that $H_{1}/H_{2}$ is cyclic of 
prime degree, and if we may continue in this way until we finally 
reach a normal subgroup $H_{n}$ which is itself cyclic of prime 
degree, we say that the group $G$ is {\em solvable}.

There are now two cases to be considered.  The group $G$ may have an 
order which is prime to the characteristic of $F$; this will always be 
the case when $F$ is of characteristic 0.  Or it may be divisible by 
the finite characteristic $p$ of $F$, in which case some of the 
composition factors will be of degree $p$.

In the first case, for each composition factor of order $q$, we may 
by definition solve the equation $x^{q}-1$ in radicals, and hence 
we may add $q$-th roots of unity whenever they are not present.  
It is a basic fact of algebra (``Kummer theory'') that in a field $L$ 
of characteristic other than $p$ containing the $q$-th roots of unity, 
any cyclic extension of degree $q$ can be given by a radical; that is, 
by the roots of an equation $x^{q}=t$, with $t \in L$.  Each of the 
composition factors can therefore be related to a corresponding cyclic 
extension defined by a radical, and in this way we can construct an 
extension containing the field $K$ derived by adjoining the roots of 
$P$ to $F$, and so factor $P$ in terms of radicals; which is to say, 
solve it in terms of radicals.  Since any polynomial which can be 
solved in radicals can be solved in Chebyshev radicals, it follows 
that any polynomial with solvable Galois group prime of order prime to 
the characteristic of the base field can be solved in Chebyshev 
radicals as well.

If the group $G$ is divisible by the characteristic $p$ of $F$, then 
$G$ will have composition factors which are cyclic of order $p$.  
These will correspond to cyclic field extensions of degree $p$.  It is 
a basic fact of algebra (``Artin-Schreier theory'') that in a field 
$L$ of characteristic $p$, any cyclic extension of degree $p$ can be 
given by a root (and hence all roots) of an equation $x^{p}-x = t$, 
with $t \in L$.  It follows that if this equation can be solved in 
radicals (ordinary or Chebyshev), then any solvable extension in 
characteristic $p$ can be solved in radicals.  

In this connection we have the following theorem.
\begin{theorem}
In a field of characteristic two, we may solve the polynomial
$x^{2}+x+t$ in Chebyshev radicals.
\end{theorem}
Proof: 
As before, we note that if $x^{2}+x+1 = 0$ then $x^{\cp 5} = 
x(x^{2}+x+1)^{2}=0$, and hence $x^{2}+x+1$ is solvable in Chebyshev 
radicals.  Hence we may assume that $t \ne 1$.

Suppose that $F_{0}$ is a field of characteristic two with $t \in 
F_{0}$.  We define $F_{1} = F_{0}(s)$, where $s^{\cp 2} = t+1$.  
Following this, we define $F_{2} = F_{1}(r)$, where $r \ne 1/s$ and 
$r^{\cp 3} = t/s^{3}$.

Now let $w=sr$. We claim that $w^{2}+w + t =0$. Since 
$r \ne 1/s$, we may multiply by $rs+1$ without introducing
any zeros. We then obtain
\[(w^{2}+w+1)(rs+1) = s^{3}r^{3}+(t+1)sr + t,\]
and on dividing by $s^{3}$ we obtain
\[r^{3}+r+t/s^{3} = r^{\cp 3}+t/s^{3} = 0.\]
$\Box$

It follows that in characteristic two, any extension $K/F$ for which
the separable closure of $F$ in $K$ has a solvable Galois group
over $F$ may be expressed in terms of Chebyshev radicals.

We also have the following negative result. 
\begin{theorem}
If $F$ is a field of characteristic $p$, where $p$ is an odd prime,
and if $t$ is transcendental over $F$, then $x^{p}-x-t$ cannot
be solved in either ordinary or Chebyshev radicals over 
$F(t)$.
\end{theorem}
Proof:
We essentially duplicate an argument already given.
Let $F$ be a field of characteristic $p$ containing the $q$-th roots 
of unity for every prime $q /ne p$, let $t$ be an element 
transcendental over $F$, and let $K_{0} = F(t)$.  

Supposing $K_{n}$ is a field obtained as a tower of $q$-th root 
extensions of prime degree over $K_{0}$ in which $x^{p}-x-t$ may be 
factored.  If $K_{n}/K_{n-1}$ is a radical extension, so that $K_{n} = 
K_{n-1}(r)$ with $r^{p} \in K_{n-1}$, then the polynomial can be 
factored over $K_{n-1}$, since $x^{p}-x-t$ is separable.  On the other 
hand, if $K_{n} = K_{n-1}(r)$ where $r^{q} \in K_{n-1}$ and $q \ne p$, 
then either $r \in K_{n-1}$ and so $K_{n} = K_{n-1}$, or $K_{n}$ is 
cyclic of degree $q$ over $K_{n-1}$, in which case $x^{p}-x-t$ must 
already be factorable over $K_{n-1}$, since there exists no 
subextension.  Iterating this process, we conclude that $x^{p}-x-t$ 
is factorable over $K_{0}$, which is a contradiction.
$\Box$

We defined solvability for a group $G$ in terms of a composition 
series; that is, a series of subgroups, each of which is normal in the 
preceding subgroup, which descends to the trivial group, and which 
has quotients which are cyclic groups of prime power order.  An 
equivalent definition relaxes the condition on the quotient to merely 
require that these be abelian, and it is often advantageous to consider 
abelian extensions (corresponding to these abelian quotients) in 
general.

Let us suppose that $P$ is a separable polynomial with coefficients in 
a field $F$ and that $L$ is a splitting field for $P$ obtained by 
adjoining all of the roots of $P$ to $F$.  Let $G$ be the Galois group 
of the extension $L/F$, $H$ a normal subgroup with an abelian quotient 
$E = G/H$, and $K$ the subfield of $L$ left fixed by $H$.  Let $n$ be 
the annihilator of $E$, that is, the least integer such that $ne = 0$ 
for every $e \in E$.  Let $n$ be prime to the characteristic of $F$, 
and let $\zeta$ be a primitive $n$-th root of unity.

In these circumstances, the Galois group of $L(\zeta)/K(\zeta)$ is 
abelian, being a quotient of $E$.  If we can solve $P$ in radicals 
over $K(\zeta)$, we will have completed the first layer of a complete 
solution in radicals.  We may begin by factoring $P$ over $K(\zeta)$; 
hence we may assume that we wish to solve an irreducible polynomial 
with abelian Galois group over a field whose characteristic does not 
divide the order of the group and which contains the $n$-th roots of 
unity, where $n$ is the annihilator of the group.  If we can do that, 
we may iterate the process until the equation is completely solved.

Let us therefore reassign our names, so that $G$ is now an abelian
group with annihilator $n$, and $F$ is a field with characteristic
prime to $n$ containing the $n$-th roots of unity. We make the
following definition.
\begin{definition}
Let $F[G]$ be the set of maps $r:F \rightarrow G$. We make this into
a vector space over $F$ by defining addition and scalar multiplication
in the obvious way, so that $(r+s)(g) = r(g)+s(g)$ for $g \in G$
and $r, s \in F[G]$; and $(\lambda r)(g) = \lambda r(g)$ for
$\lambda \in F$. We then define a product on $F[G]$ by setting
\[(r \ast s )(g) = \sum_{e+f=g}r(e)s(f).\]
The $F$-algebra so defined is called the {\em group algebra}
of $G$ over $F$.
\end{definition}

We may define another product on $F[G]$, making it an algebra
in another way, by means of coordinatewise multiplication.
\begin{definition}
For $r, s \in F[G]$, we define
\[(r \cdot s)(g) = r(g)s(g).\]
\end{definition} 

Of particular interest are the elements of $F[G]$ which are
homomorphisms of $G$ to the multiplicative group $F^{\times}$
of $F$.
\begin{definition}
An element $\chi \in F[G]$ is called a {\em character} of $G$ if
$\chi (g+h) = \chi(g) \cdot \chi (h)$.
\end{definition}

The characters of $G$ together with the product ``$\cdot$'' form an 
abelian group denoted by $\hat G$. We have a natural and 
nondegenerate pairing between $G$ and $\hat G$ given by
\[ (g, \chi) \mapsto \chi(g).\]
Using this, we may conclude that $G$ is isomorphic to $\hat G$,
and canonically isomorphic to $\hat {\hat G}$. This is therefore
a duality, and $\hat G$ is called the {\em dual group} to $G$.

To understand the structure of a commutative algebra, it is helpful
to identify idempotent elements. In the case of $F[G]$ we have
the following.

\section{Class Field Theory}

\subsection{The $P$-adic Chebyshev Radical}

We have seen the great utility of the Chebyshev exponent as a means
of expressing the roots of solvable polynomials. The solution of 
polynomial equations in Chebyshev exponents is interesting in another
way: it sometimes happens that arithmetical questions are more easily 
analysed in Chebyshev terms. In particular, it is helpful when 
analyzing the unramified extensions which class field theory associates
with the ideal class group.

When studying arithmetical questions, it is often of great benefit to 
work locally, over the completions of the number field.  We will 
recall for the reader of some of the basic facts about $p$-adic 
fields; the reader who requires more should consult the standard 
references, such as \cite{Kob}.

We make the following standard definition.
\begin{definition}
Let $p$ be a prime number, and let ``$| |_{p}$'' denote a function from
$\Q$ to $\R$ such that $|x|_{p} = 0$ if and only if $x = 0$, and
such that if $q = p^{n}r$ is a rational number, where $r$ is a rational
number such that neither the numerator nor the denominator is divisible
by $p$, then $|q|_{p} = p^{-n}$. Then the $p$-adic completion of $\Q$,
$\Q_{p}$, is the topological field $C_{p}/N_{p}$, where $C_{p}$
is the ring of Cauchy sequences for the metric defined by
$|x - y|_{p}$, and $N_{p}$ is the maximal ideal of sequences 
converging to 0 under this metric.
\end{definition}

By the way this is defined, all Cauchy sequences converge, and hence 
$\Q_{p}$ is complete.  It is not hard to see that $N_{p}$ is maximal, 
and hence that $\Q_{p}$ is indeed a field.  It differs from the field 
of real numbers in that the algebraic closure of $\Q_{p}$ is not 
topologically complete.  To find a $p$-adic field which plays the same 
role as the complex numbers do, it is usual to complete the topological 
closure, obtaining thereby a field both algebraically closed and 
topologically complete.
\begin{definition}
We denote by $\C_{p}$ be the topological completion of the
algebraic closure of the $p$-adic field $\Q_{p}$. The function
from $\C_{p}$ to the reals extending the absolute value on
$\Q_{p}$ to $C_{p}$ we will again denote by $| |$; we retain
the normalizing condition $|p|_{p} = \frac{1}{p}$.
\end{definition}

We began by defining Chebyshev powers over arbitrary rings, so
of course they are defined over any $p$-adic field. By using
the same power series we derived over $\C$, we may also define
Chebyshev functions for $p$-adic values of the exponent.

\begin{definition}
For any $x$ and $k$ in $\C_{p}$, we define 
$$(2+x)^{\cp} = 2+k^{2}x+\frac{2}{4!}k^{2}(k^{2}-1)x^{2}+
\frac{2}{6!}k^{2}(k^{2}-1)(k^{2}-4)x^{3}+ \cdots,$$
whenever the series above, which is (\ref{pows2}), converges.
\end{definition}

\begin{theorem}
Let $F$ be a topologically complete subfield of $\C_{p}$, and let
$k$ be an element of $F$ such that$|k|_{p}>1$. Then for any $x \in F$
the power series (\ref{pows2}) defines an element $x^{\cp k} \in F$
if and only if
\begin{equation}
\label{rad1} 
|x-2|_p < 1/(|k|_p \sqrt[p-1]{p})^{2}.
\end{equation}
\end{theorem}
Proof:
If $|k|_{p}>1$ then $|k+i|_{p} = |k|_{p}$ for any integer $i$.
Hence (\ref{pows2}) converges precisely when
$$2 + \frac{2}{2!}k^{2}x + \frac{2}{4!}k^{2}x^{2}+ \cdots$$
converges. However, this is the power series expansion of
$2 \cosh (k \sqrt{x})$. Since (\cite{Kob}) $\cosh x$ converges
if and only if 
$$|x|_{p} < 1 / \sqrt[p-1]{p},$$
(\ref{pows2}) will converge if and only if
$$|k \sqrt{x}|_{p} < 1 / \sqrt[p-1]{p},$$
which proves (\ref{rad1}).
$\Box$

\begin{theorem}
Let $F$ be as before. Then for
any $x \in F$ and $k \in \Z_{p}$,
(\ref{pows2}) converges and defines
$x^{\cp k} \in F$ when 
$$|x-2|_{p} < 1.$$
\end{theorem}
Proof:
Since (\ref{pows2}) gives a polynomial with integer coefficients
when $k$ is a positive integer, the coefficients of (\ref{pows2}) are
integers for positive integer values of $k$. Hence they are less than or
equal to one in absolute value if $k \in \Z_{p}$ by the density of
the positive integers in $\Z_{p}$. Hence if $|x|_{p} < 1$ the terms
of (\ref{pows2}) form a null sequence, and so the series is convergent.
$\Box$

\begin{theorem}
Let $F$ be as before, and let $x$, $n$, and $m$ be elements of
$F$. If $(x^{\cp n})^{\cp m}$ is defined, then $x^{\cp nm}$ is defined
and $(x^{\cp n})^{\cp m} = x^{\cp nm}$
\end{theorem}
Since this is an identity over $\C$ when $|x-2| < 2$, it is an identity
in formal power series. Hence when the series converge, it must remain
true in $F$.
$\Box$

The corresponding results for $\Sc_{k}(x)$ follow immediately, and
may easily prove corresponding results for $\Uu_{k}(x)$.

\begin{definition}
For $x$ and $k$ in $\C_{p}$, we define
$$\Uu_{k}(2+x) = = k + \frac{k(k^{2}-1^{2})}{3!}(\frac{x}{4})+
\frac{k(k^{2}-1^{2})(k^{2}-3^{2})}{5!}(\frac{x}{4})^{2}+\cdots,$$
whenever the power series above, which is (\ref{eq:u-series}), converges.
\end{definition}

\begin{theorem}
Let $F$ be as before, and let $k$ be an element of $F$ such that
$|k|_{p}>1$. Then $\Uu_{k}(x)$ is defined if and only if
\begin{equation}
\label{rad2} 
|x-2|_p < |4|_{p}/(|k|_p \sqrt[p-1]{p})^{2}.
\end{equation}
\end{theorem}
Proof:
We argue as before; the series converges precisely when the
series for $\sinh (k \sqrt{x/4}) / \sqrt{x/4}$ converges,
which entails
$$|k \sqrt{x/4}|_{p} < 1 / \sqrt[p-1]{p},$$
and hence (\ref{rad2}).
$\Box$

\begin{theorem}
Let $F$ be as before. Then for
any $x \in F$ and $k \in \Z_{p}$,
(\ref{eq:u-series}) converges and defines
$\Uu_{k}(x) \in F$ when 
$$|x-2|_{p} < |4|_{p}.$$
\end{theorem}
Proof:
We argue as before, except now we note that
(\ref{eq:u-series}) gives a polynomial with integer coefficients
when $k$ is an odd positive integer, and odd positive integers
are dense in $\Z_{p}$ when $p$ is odd. When $p=2$, the odd positive
integers are dense in the odd 2-adic integers; since these are larger
than the even 2-adic integers the result follows in this case also.
$\Box$

\subsection{Showing an Extension is Unramified}

We are interested in finding unramified abelian extensions.  One
means for accomplishing this is to use of pure
radical extensions.  In this instance, we have the following well-known
criterion:  a pure radical extension $F(\sqrt[n]{a})$
of a number field $F$ is unramified outside of $n \infty$ if and
only if the principal ideal $(a)$ generated by $a$ is an $n$-th
power as an ideal.  In particular, this includes extensions generated
by roots of units of $F$.  This is easily seen to be the case by
working locally, where $a$ becomes an $n$-th power times a unit.

We will find it convenient to work with Chebyshev radical extensions
as well.  Similar conditions to the above in terms of the trace of an
ideal from a certain field can be developed.

We will be particularly interested in split unramified extensions.
Suppose we have a Galois extension $L/F$ of a number field $F$
whose Galois group is split with an abelian kernel $A$.
We then have a subextension $K/F$ with Galois group $G$, and
the Galois group of $L/F$ is
$\Gal(L/F) \iso G \rtimes A$ of $L/F$.
The extension $L/F$ can be given as the splitting field of
a polynomial $g$ of degree $|A|$ with coefficients in $F$, 
and $K/F$ as the splitting field of a polynomial $f$ of lesser degree.

In such a situation,
the extension $L/K$ is unramified at primes
lying over the prime $\wp$ of $F$ if
and only if the inertia groups for the commutative algebra given by the
roots of $f$ over $F_{\wp}$ are the same as those for the roots of $g$
over $F_{\wp}$.  

$\Gal (K/F)$ as a permutation group on the roots of $f$
injects into $\Gal (L/F)$ as a permutation group on the roots of $g$.
If $r_f$ is a root of $f$ and $r_g$ is a root of $g$,
for $\wp$ a prime in $F$ which ramifies totally in $F(r_f)$, 
we have a correspondingly that the factorization of
$\wp$ in $F(r_g)$ will consist of unramified primes, and
$m$ primes which ramify with inertia corresponding to 
the inertia of $\wp$ in $F(r_f)$.  Hence locally at
$\wp$, we have that the $\wp$ factor in the discriminant of
the field $F(r_g)$ is the $m$-th power of what it is
for $F(r_f)$.

\subsection{Unramified Families For Degree Three}

The simplest example of the situation described above is where $G$ is
cyclic of degree two, and $A$ cyclic of degree three.  In this case we
are seeking extensions of $L/F$ with Galois group $S_3$, such that if
$K$ is the quadratic subfield, then $L/K$ is unramified.  Looking at
degree three polynomials, this means that at all primes $\wp$ of $F$
which ramify in $L$ we do not get total ramification in the degree
three field $F(r_g)$, but rather a product of a prime of ramification
index one and one of index two, the latter corresponding to the (total)
ramification of $\wp$ in $K$.  Hence, $L/K$ will be unramified if and
only if the field discriminants of $F(r_f)$ and of $F(r_g)$ are the
same.  This happens if and only if no prime ramifies with index three
in $F(r_g)$.

The polynomial
$$x^3+bx+c$$
is {\em generic} for the Galois group $S_3$. By this I mean that
over $\Q (a,b)$ its splitting field has group $S_3$, and 
(in characteristic other than 3) any Galois extension with
group $S_3$ can be obtained as the splitting field of a polynomial
of this form.

Let us therefore consider a polynomial
$$x^3 + bx + c,$$
with coefficients in a number field $F$.

If we set $z = \sqrt{-3/b}x$, then we find that $z$
satisfies the polynomial
$$z^3 - 3z + \sqrt{{-27} \over b^3}c,$$
so that if $h = - \sqrt{{-27}/b^3}c$, we have a solution in
Chebyshev radicals, given by
$$z_1 = \sqrt[\cp 3]{h},$$
$$z_2 = -\sqrt[\cp 3]{-h},$$
and
$$z_3 = -z_1 -z_2.$$

In order not to trouble ourselves with square roots
we can eliminate them by using instead the polynomial for
$$z^{\cp 2} = -{{3x^2} \over b}-2.$$
Up to a change of scale, this transformation is what
results from squaring the Lagrange resolvents and finding roots
in terms of these.

If the polynomial is irreducible, then the splitting field for
$z^{\cp 2}$ gives the same extension as for $x$.  On the other hand,
$$z^{\cp 2} = (\sqrt[\cp 3]{h})^{\cp 2} = \sqrt[\cp 3]{h^{\cp 2}} =$$
$$\sqrt[\cp 3]{-2- {{27c^2} \over b^3}}
= \sqrt[\cp 3]{2+{{27c^2} \over b^3}},$$
which gives us a root when it is defined.

If we denote the discriminant of the original polynomial in
$x$ by 
$$\d = -4b^3 -27c^2,$$
then we can write
$$-2 - {{27c^2} \over b^3}$$
as 
$$2 + {\d \over b^3},$$
so that we have a root given by
$$\sqrt[\cp 3]{2 + {\d \over b^3}},$$
when this is defined.

If $\wp$ is a prime in a number field $F$ lying over the
rational integral prime $p$, we normalize
the associated valuation in the usual way, by setting
$$|p|_{\wp} = |p|_p = 1/p.$$
We now can state a criterion for determining unramified
extensions.

\begin{theorem}
If $x^3 + bx +c$ with $b \ne 0$ has 
coefficients in an algebraic number field $F$ and is
irreducible over $F$,
and if $\e = -2 - 27c^2/b^3$, $\d = -4b^3-27c^2$, and $r$ 
is any root,
then $F(r, \sqrt{\d})$ gives
an unramified extension of $F(\sqrt{\d})$ if one of the following
three conditions is true for each prime $\wp$ which ramifies
in $F(\sqrt{\d})$, and gives an
extension unramified outside of primes lying over 3 if and only if
one of the three conditions always holds.  
If the field $F$ contains no subextension which is
unramified over $\Q$ at 3, then the
conditions are both necessary and sufficient.

\begin{enumerate}

\item

$$|{{\e+2} \over {27}}|_{\wp} = |{c^2 \over b^3}|_{\wp} < 1,$$

\item

$$|{{\e-2} \over {27}}|_{\wp} = |{\d \over {27 b^3}}|_{\wp} < 1,$$

\item

$$|{{\e^2} \over {27}}|_{\wp} = |{{(\d + 2b^3)^2} \over {27b^6}}|_{\wp} < 1.$$
\end{enumerate}
\end{theorem}
Proof:
Since we have complex embeddings of $F(x)$ only if we have
corresponding complex embeddings of $F(\sqrt \d)$, we never 
have ramification at an infinite place.
Hence we will have an unramified extension if and only if
there is a $\wp$-adic root of the above polynomial for
each prime $\wp$ of $F$ which ramifies in $F(\sqrt \d)$.
And we will have such a root if and only if we have a root of 
$x^2 -3x - \e$, where $\e = -2 - 27c^2/b^3$.

If $\wp$ is prime to 3, then
$-T_{1 \over 3}(-\e)$ will give us a root if
$$|\e+2|_{\wp} = |{{27c^2} \over b^3}|_{\wp} = {|{c^2 \over b^3}|_{\wp}} < 1$$
by the lemma on the power series expansion of $T_n$.
On the other hand, if $\wp$ lies over 3, then we will
have a root if
$$|\e+2|_{\wp} = |{{27c^2} \over {b^3}}|_{\wp} < 1/27,$$
which entails
$$|{c^2 \over b^3}|_{\wp} < 1.$$
So if the first condition is true for a prime $\wp$, then
we do not have ramification at $\wp$.

Since 
$$\e = -2-{{27c^2} \over b^3} = 2 + {\d \over b^3},$$
we will have a root $\sqrt[\cp 3]{\e}$ when

$$|{\d \over b^3}|_{\wp} = |{{27 \d} \over b^3}|_{\wp} < 1$$
for $\wp$ prime to $p$, and
$$|{\d \over b^3}|_{\wp} < 1/27.$$
for $\wp$ lying over $p$.  Putting these together, we
get the condition
$$|{\d \over {27b^3}}|_{\wp} < 1.$$
So if the second condition is true for a prime $\wp$, we
do not have ramification at $\wp$.

Finally, if $|\e|_{\wp} < 1$ and $\wp$ does not lie over 3, then
$\sqrt[\cp 3]{2-\e^2}$ will give a root of the polynomial
transformed for the second time by $x \mapsto x^{\cp 2}$.  And if
$|\e^2|_{\wp} < {1 \over {27}}$ for a prime $\wp$ lying over 3, then
$T_{1 \over 3}(2-\e^2)$ will likewise give a root.  Putting these
conditions together, we have a root of the transformed polynomial
when $|{\e^2}/{27}|_{\wp} < 1$, so the third condition
entails we do not have ramification at $\wp$.

If one of these conditions is true for each prime
$\wp$ which ramifies in $F(\d)$, then we have an unramified
extension at every prime, and hence a globally unramified
extension.  This is the first part of the theorem.

For the second part of the theorem, it is helpful to
assume that the polynomial is $\wp$-reduced.
The conditions are in terms of $c^2/b^3$ and $\d/b^3$, and
so remain the same if we make the transformation
$b \mapsto k^2b$, $c \mapsto k^3c$.  Hence we may assume the
polynomial is $\wp$-reduced, meaning that $|b|_{\wp} \le 1$,
$|c|_{\wp} \le 1$, and either $|b|_{\wp} > |\wp|_{\wp}^2$, or
$|c|_{\wp} > |\wp|_{\wp}^3$, or both.

Assuming this reduction, we can have $b$ a unit and
$c$ a non-unit, $b$ and $c$ both units, or $c$ a unit
and $b$ not a unit.

If $b$ is a unit and $c$ is not, then the first condition holds.

If both $b$ and $c$ are units, then since $|\d|_{\wp} < 1$
and $p \not\mid 3$, the second condition holds.

If $c$ is a unit and $b$ is not, then since $\wp$ divides 
$\d =-4b^3-27c^2$, $\wp$ must divide $27c^2$, and since $c$ is a unit,
$\wp$ must lie over 3.  If $F/\Q$ contains no subextension which is
unramified at 3, then if $r$ is a root of the $\wp$-reduced polynomial
$x^3+bx+c$ which does not lie in $\Z_p$, then one of the conditions
$|r+2|_{\wp} < {1 \over 3}$, $|r|_{\wp} < {1 \over 3}$, or
$|r-2|_{\wp} < {1 \over 3}$ must hold.  Since $\e = r^{\cp 3}$, then this
entails one of the conditions $|\e+2|_{\wp} < {1 \over {27}}$,
$|\e|_{\wp} < {1 \over 9}$, or $|\e-2|_{\wp} < {1 \over {27}}$
must hold; which then shows that one of $|(\e+2)/27|_{\wp}< 1$,
$|\e^2/81|_{\wp} < 1$ or $|(\e-2)/27|_{\wp} < 1$ must hold.  
Hence one of the conditions holds, and this is the
second part of the theorem.
$\Box$

We may now employ this theorem to construct families of unramified
extensions.

\begin{corollary}
Let $$x^3 + bx + b^2t$$ be an irreducible polynomial
with $b$ and $t$ in the ring of integers of a number
field $F$. Let $r$ be a root and 
let $d=-27bt^2-4$.  Then $F(r, \sqrt{bd})$ is an
unramified extension of $F(\sqrt{bd})$.
\end{corollary}
Proof:
The first condition holds if $|tb^2|_{\wp} < 1$, and hence holds if
$b$ is not a $\wp$-unit.  If $b$ is a $\wp$-unit, then 
$|d|_{\wp} < 1$, since $\d = b^3d$.  Hence if $\wp$ is prime to 3,
$$|{\d \over {27b^3}}|_{\wp} = |{d \over {27}}|_{\wp} < 1.$$
But we cannot have $\wp$ lying over 3, 
since $\d \equiv -4b^3 \pmod {27}$
and $b$ is a $\wp$-unit, and so $\d$ is prime to 3.
Hence the second condition holds if the first condition fails,
and we have an unramified extension in all cases.
$\Box$

\begin{theorem}
Let $F$ be a number field, and $b$ an non-zero integer in $F$.
Then
$$x^3 + bx + c$$
is a congruence family.
\end{theorem}
Proof:

Let $c$ be an integer of $F$ such that $x^3+bx+c$ produces an
unramified extension.
By the previous corollary, values of $c$ of the form
$tb^2$ will produce unramified extensions, so such $c$ exist.

Now let $h$ be an integer of $F$ such that 
$$|{h \over {27b^3}}|_{\wp} < 1$$
for every prime $\wp$ of $F$ dividing $(3b)$.
Let $\delta = -4b^3-27c^2$ be the discriminant of the polynomial
$x^3+bx+c$, and $\delta' = -4b^3-27(c+h)^2$ the discriminant
of $x^3+bx+c+h$.  For any prime $\wp$ dividing $\delta'$, we will
have a local root if 
$$|{{\delta'} \over {27b^3}}|_{\wp} < 1,$$
so we may confine our attention to primes where
$$|{{\delta} \over {27b^3}}|_{\wp} \ge 1.$$

In this case, since
$$|{{\delta-\delta'} \over {27b^3}}|_{\wp}=|{{h(h+2c)} \over {b^3}}|_{\wp} < 1,$$
we have also
$$|{{\delta} \over {27b^3}}|_{\wp} \ge 1,$$
and hence $|h|_{\wp} < |\delta|_{\wp}$.

The $\wp$-adic root $r$ of $x^3+bx+c$ will converge by 
Newton's method to a $\wp$--adic root of $f(x) = x^3+bx+c+h$ if
$$|f(r)|_{\wp} < |f'(r)|_{\wp}^2,$$
which means precisely if
$$|h|_{wp} < |\delta|_{\wp},$$
and hence we get a $\wp$-adic root in all cases, and
hence an unramified extension.
$\Box$

While eg. the modulus $81b^4$ will always work, it is
by no means the case that so large a large modulus is always needed.  For
instance, we have the following, when the base field $F$ is $\Q$:

\begin{corollary}
Let $b$ be $\pm p$, for $p$ a rational integral prime different from 3.
Then $x^3+bx+c$ for integral $c$ gives 
an unramified extension of its quadratic
subfield if it is irreducible, and if $c$ is either prime to $p$, or
divisible by $p^2$.
\end{corollary}
Proof:
Let $\ell$ be a prime dividing $\delta = -4b^3-27c^2$.  Then 
since $\delta \equiv -4b^3 \pmod 3$, $\ell \ne 3$.  If $\ell \ne p$,
we have
$$|{{\delta} \over {27b^3}}|_{\ell} < 1.$$
Hence we will only possibly have ramification if $\ell = p$.
If $c$ is divisible once by $p$, then we have total ramification
at $p$.  If $c$ is divisible more than once, we have
$$|{{c^2} \over {b^3}}|_p < 1,$$
and hence no ramification.
$\Box$

In this case, we see that the minimum modulus defining the
congruence family is $p^2$.  If $p=3$, we can do a similar analysis;
however as we shall see there is another way of looking at this
case which has an interesting generalization.

\begin{theorem}
Let $b$ be a fixed nonzero rational integer.  Then the polynomial 
$$x^3 + bx + c$$
for rational integer values of $c$ such that the above polynomial is
irreducible will give an unramified extension of $\Q(\sqrt{\d})$ 
if and only if mod $b^3$ we have 
$$c^2 \equiv 0 \pmod {b^3}$$
or 
$$c^2 \equiv -{{4b^3} \over {27}} \pmod {b^3}.$$
\end{theorem}

\subsection{Unramified Families for Degree Four}

Just as 
$$x^3+bx+c$$ is generic in characteristic other than 3 for
the triangle group $D_3 = S_3$, the polynomial
$$x^4+bx^2+c$$
is generic for the square group $D_4$ in characteristic other
than 2.  It has a splitting field with Galois group $D_4$
over $\Q (a,b)$, and outside of characteristic 2, any Galois
extension with group $D_4$ can be obtained as the splitting field
of a polynomial in this form.

The reduction to this form is a little less trivial.  If we have
a polynomial of degree four with roots $r_0,r_1,r_2,r_3$ giving
a Galois extension with group $D_4$, we may obtain another polynomial
with roots $r_0-r_2$ and its conjugates, if the roots
are correctly ordered.  One way to obtain this polynomial is
to factor a resolvent polynomial of degree twelve.

\begin{lemma}
If 
$$x^4+a_1x^3+a_2x^2+a_3x+a_4$$
is a polynomial with coefficients in a field $F$ of characteristic
other than 2, whose splitting field gives a Galois extension
of $F$ with group $D_4$, then
$$z^{12}-3 a1^2z^10+8 a2z^{10}$$
$$+(-16 a1^2 a2+8 a4+3 a1^4-2 a3 a1+22 a2^2)z^8$$
$$+(26 a3^2+8 a1^4 a2+16 a4 a2- a1^6-24 a1^2 a2^2-30 a3 a1 a2$$
$$+28 a2^3+8 a3 a1^3-6 a4 a1^2)z^6$$
$$+(17 a2^4+48 a3^2 a2-25 a3^2 a1^2+6 a4 a1^4-54 a3 a2^2 a1$$
$$-112 a4^2-12 a1^2 a2^3+38 a3 a1^3 a2-6 a3 a1^5+24 a2^2 a4$$
$$-32 a4 a1^2 a2+2 a1^4 a2^2+56 a3 a1 a4)z^4$$
$$+(6 a3 a1^3 a2^2+216 a3^2 a4-120 a3 a1 a4 a2+72 a4^2 a1^2+42 a3^2a1^2a2$$
$$+18 a3^2 a2^2+18 a3 a1^3 a4-9 a3^2 a1^4-54 a3^3 a1+32 a4 a2^3$$
$$-26 a3 a1 a2^3- a2^4 a1^2z^2-192 a4^2 a2+4a2^5z^2-6 a4 a2^2a1^2)z^2$$
$$-27 a1^4 a4^2-192 a3 a1 a4^2-4 a3^3 a1^3-80 a3 a1 a2^2 a4-128 a2^2a4^2$$
$$+ a2^2 a3^2 a1^2+256 a4^3-4 a2^3 a1^2 a4-27 a3^4+16 a2^4 a4-6 a3^2a1^2 a4$$
$$-4 a2^3a3^2+144 a2 a1^2 a4^2+144 a4 a3^2 a2+18 a3 a1^3 a2 a4+18a3^3a1 a2$$
will factor into either a single factor of degree four and a single
factor of degree eight, or into a single factor of degree four
and a repeated factor of degree four.  The single factor of degree
four is of the form
$$z^4+bz^2+c$$
and is such that if $r$ is a root of our original polynomial, then
there is a root $s$ of this polynomial such that $F(s)$ defines the
same field extension as $F(r)$.  
\end{lemma}
Proof:

If the Galois group is $D_4$, it will contain a generator
$\sigma$ for the cyclic subgroup of order four.  Applying this to a root
$r_0$, we obtain $r_i = r_0^{\sigma^i}$.  Then $s_0=r_0-r_2$ and
its conjugates generate the same extension, and it is trivial 
to check that the polynomial derived from these roots has a
zero $x$ and $x^3$ term.  Moreover, $F(s_0)$ generates the same
extension as $F(r_0)$, since this is the same extension as
$F(r_2)$.

The polynomial of degree twelve in the statement is the 
polynomial for $r_i - r_j$; it is 
readily calculated using e.g. resultants.  Hence it will have
a factor of degree four corresponding to the polynomial considered
in the previous paragraph. 
$\Box$

Let $\d = b^2-4c$.  In some cases, the field $L=F(r_0,r_1,r_2,r_3)$  
obtained by adjoining the roots of 
$$x^4+bx^2+c$$
will define an unramified extension of
$K=F(\sqrt{c\d})$, where $b,c \in F$. We consider first
infinite places.

\begin{theorem}
The extension $L/K$ is unramified at any real infinite place 
iff either 
$b<0$ and $c>0$, or $\d < 0$ at that place.
\end{theorem}
Proof:

If $\d < 0$ at the place in question, $K$ is already complex
there and the extension is unramified trivially.  The polynomial
$x^2+bx+c$ 
has two positive real roots if and only if $b<0$ and $c>0$, and
in this case, $L$ is real at this place and hence the extension
is unramified there.
$\Box$

\begin{theorem}
If $x^4+bx^2+c$ is a polynomial with coefficients in a number
field $F$ whose roots generate a $D_4$ extension of $F$, and
if $r$ is any root and $\d = b^2-4c$, then $L=K(r)$
is unramified at a prime $\wp$ of $K=F(\sqrt{c\d})$ lying over
a prime $p$ of $F$ if one
of the following conditions holds:

\begin{enumerate}

\item

The prime $p$ divides $c\d$ an odd
number of times, and
$$|\frac{c}{b^2}|_p < 1.$$

\item

The prime $p$ divides $c\d$ an odd number of times, and
$$|\frac{\d}{2b^2}|_p < 1.$$

\end{enumerate}

\end{theorem}
Proof:

If $|x| < 1$ in any valuation, we have that
$$\sqrt{1-4x} = 1-2x - 2x^2 \cdots$$
converges.  Hence if $|\frac{c}{b^2}|_p < 1$, we have that
$$\sqrt{b^2-4c} = b \sqrt{1-4\frac{c}{b^2}}$$
has a $p$-adic root, and hence the roots of 
$$x^4+bx^2+c$$
are either elements of $F_p$ of square roots of such elements.  The
ramification is thus locally no more than degree two.  If $p$ divides
$c\d$ an odd number of times, the ramification of $F(\sqrt{c\d})$ is
also two, and hence no further ramification can occur.

If we set $s_0 = r_0 - r_1$ together with its conjugates, we obtain
the four roots of the polynomial
$$x^4 + 2bx^2 + \d.$$
Hence the second part of the theorem follows from applying the conclusions
of the first part to this polynomial, which defines different stem fields
but the same splitting field.
$\Box$

\begin{corollary}
Let
$$x^4+bx^2+b^3t$$
or
$$x^4-bx^2+b^3t$$
be a polynomial with $b \ne 0$ and $b$ and $t$ 
in the ring of integers of a number field $F$, defining a $D_4$ 
extension of $F$.
Let $r$ be a root 
and let $\d = 1-4bt$.  Then
$F(r, \sqrt{b\d})$ is an unramified extension of $F(\sqrt{b\d})$
at all non-infinite places.
\end{corollary}
Proof:

We have
$$(\pm b)^2 - 4b^3t = b^2\d.$$
If $p$ is a finite prime of $F$ which divides $b$, then
it does not divide $\d$.
$$|\frac{b^3t}{b^2}|_p = |bt|_p < 1,$$
and hence 
$S_3$ extensions of any number field $F$ (this is obvious both from
its derivation and from the fact that specializing to $u=1$ leads us
to a polynomial which is still generic.)  And it is well-adapted to
the purpose of analyzing when a degree extension is unramified over
the quadratic subfield, which is given by the roots of
$x^2+u(4s^3+27t^2u)$.  We may analyze these by fixing $s$, and we begin
by setting $s=1$.

In this case, the basis $[1,x,x^2/u]$ consists of elements integral over
$\Z[u,t]$, and the trace form of this basis has discriminant
$-u(27ut^2+4)$.  The polynomial giving the quadratic subfield
is $x^2 + u(27ut^2+4)$.  If $u$ is even, then $[1,x/2]$ removes
a factor of 4, whereas if $u$ is odd, $[1,(x+t)/2]$ removes
a factor of 4, so we are
left with $-u(27ct^2+4)$ in all cases.  

If this is square free we
are done, and we have the same discriminant for the degree three
field as for the degree two field.  If the discriminant contains a 
factor which is the square of a prime $p$ which is prime to $2u$,
we remove may remove the square term from both the degree two
and degree field discriminants, by $[1,x/p]$ for the degree two
field, and $[1,x,(x_1-x_2)^2/p]$ (where $x_1$ and $x_2$ are the
two conjugates of $x$) for the degree three case.  In case $p$ 
divides $2u$ we modify this a bit, but in any case we simply
remove the prime factor from the degree two and degree three
extension in just the same way, which we can do since the 
ramified primes and the trace forms correspond.

Since the field will be totally real if the quadratic subfield
is real, we need not concern ourselves with the place at infinity;
and for this particular polynomial we are dealing with negative
discriminants in any case.
Hence we may conclude that substituting integer values for $t$ and $u$ into
$x^3 + ux + tu^2$ gives an unramified extension of the quadratic
subfield given by $x^2 + u(27ut^2+4)$, in case both polynomials
are irreducible.

While this result is very satisfactory in that it constructs
an infinity of unramified extensions, it is limited to doing so
over imaginary quadratic fields.  For this reason it is worth our
while to investigate other values for $s$.  Since substituting
$u \mapsto su, t \mapsto st, x \mapsto sx$ into $x^3 + ux + tu^2$
leads to $x^3 + sux + tu^2$, we can conclude that this is unramified
outside the primes dividing $s$ (where $s$ is a local unit)
by our previous analysis.  Hence all we need to do is consider 
values of $t$ and $u$ $s$-adically; and we may began by noting
that the above substitution lets us conclude immediately that
if both $u$ and $t$ are divisible by $s$, we have an unramified
extension.

If $s=2$, we obtain $x^3 + 2ux + tu^2$, with quadratic subfield
$x^2+u(27t^2u+32)$.  We know that if both $u$ and $t$ are even the
extension is unramified.  If $u$ is even and $t$ is odd, then we have
total ramification at 2 unless $u$ is divisible by $8$.  However, 
$u \mapsto 8u, x \mapsto 4x$ leads us back to $x^2 + ux + tu^2$, and
so in this case also the extension is unramified.

If $u$ is odd, then modulo 4 the polynomial becomes
$x^3 + 2x + t$.  This is totally ramified at 2 if $t \equiv 2 \pmod 4$,
and is unramified at 2 if $t$ is odd.  If $t \equiv 0 \pmod 4$, then
it is ramified at 2 but not of degree three, and so the ramification
corresponds to the ramification of the quadratic subfield.

The conclusion is that we have an unramified extension
if and only if $u$ is divisible by 8, or if it is even and
$t$ is even, or if it is odd and $t$ is not congruent to 2
mod 4.

If $s=3$, then we will obtain an unramified extension if
$u$ is divisible by three and 

If $u$ is not divisible by three,
we will obtain an unramified extension if 
$t \equiv 0 \pmod 9$; or if $u^9 \equiv -1 \pmod {27}$ and 
$ut^2 \equiv -4 \pmod {27}$; or finally if 
$u^9 \equiv 1$ and $ut^2 \equiv -2 \pmod 9$.

If $c \equiv \pm 3 \pmod 9$, then we get a polynomial with 
integer coefficients if $t = t'/9$, where $t'$ is integral.
Under these conditions, we obtain an unramified extension
if $t' \equiv 0 \pmod {27}$, or if $t' \equiv \pm c/3 \pmod 9$.

Finally, if $c \equiv 0 \pmod 9$, then transforming by
$c = 9c'$, $t = t'/3$, and $x = 3 x'$ leads to the same polynomial,
so we have already done this case.

In conclusion, we get the following

\begin{theorem}
There exist an infinite number of real quadratic fields as
well as an infinite number of imaginary quadratic fields with
3-cycles in the nongenus class group.
\end{theorem}

To find families of unramified extensions corresponding to
4-cycles in the class group of quadratic extensions, we need to
proceed a little differently, since a degree four dihedral 
extension has three quadratic subfield rather than one.
We could proceed as before by analyzing a generic polynomial
for dihedral extensions (such as $x^4-4acx^2-a^2c(b^2-4c)$
for instance), but we will content ourselves with something
simpler.

Consider the polynomial 
$$x^4 - x^3 - tx^2 - x + 1.$$
Over $\Q (t)$ this has dihedral Galois group, and three quadratic
subfields given by $x^2 - t(t-4)$, $x^2 - 4t - 9$, and
$x^2 - t(t-4)(4t+9)$.  The discriminant
of the polynomial is $t(t-4)(4t+9)^2$. 

The for any integral specialization of $t$, the ideal
$(t,t-4)$ divides 4, the ideal $(t,4t+9)$ divides 9, and
the ideal $(t-4,4t+9)$ divides 25.  Hence for any prime $p>5$
we have that it divides only one of these
three.  We can write our polynomial in three different ways:

$$(x^2+x+1)(x-1)^2 - tx^2,$$
$$(x^2-3x+1+1)(x+1)^2 - (t-4)x^2,$$
$$(x^2-x/2+1)^2 - (t+4/9)x^2.$$

This corresponds to ramification of order two and with multiplicity one
if $p$ divides $t$ or $t-4$ an odd number of times, and with
multiplicity two if $t$ divides $4t+9$ an odd number of times.

If $t$ is odd then $t(t-4)(4t+9)$ is odd.  If $t \equiv -1 \pmod 3$
then $t(t-4)(4t+9)$ is not divisible by three.  And if $t$ is
congruent to $1$ or $\pm 2$ mod 5, then $t(t-4)(4t+9)$ is not
divisible by five.  Hence for the combined congruence condition
mod 30, if $t$ is of the form $30n-13$, $30n-7$ or $30n+11$
we have that $t(t-4)(4t+9)$ is not divisible by 30, and
hence that the degree four polynomial for this specialization
is unramified over the quadratic field $\Q(\sqrt{t(t-4)(4t+9)})$

Hence we have

\begin{theorem}
There exist an infinite number of quadratic fields with both
real and imaginary discriminant and with a 4-cycle in the
class group, such that the corresponding portion of the
Hilbert class field is generated by
$$x^4-x^3-tx^2-x+1$$
for values of $t$ such that $t \equiv -13,-7,11 \pmod {30}$.
\end{theorem}

\end{document}